\documentclass[aop,MSNbibl,seceqn,dvips]{arximspdf}
\usepackage{mathbh}
\usepackage{graphicx}

\doi{10.1214/10-AOP557}
\volume{39}
\issue{2}
\pubyear{2011}
\firstpage{507}
\lastpage{548}

\makeatletter

\newtheorem{theorem}{Theorem}[section]
\newtheorem{corollary}[theorem]{Corollary}
\newtheorem{lemma}[theorem]{Lemma}
\newtheorem{proposition}[theorem]{Proposition}
\newproclaim{Remark}{Remark}

\makeatother

\begin{document}
\begin{frontmatter}

\title{Asymptotic behavior of the gyration radius for long-range
self-avoiding walk and long-range oriented percolation}
\runtitle{Asymptotic behavior of the gyration radius}

\begin{aug}
\author[A]{\fnms{Lung-Chi}
\snm{Chen}\thanksref{t1}\ead[label=e1]{lcchen@math.fju.edu.tw}}
\and
\author[B]{\fnms{Akira}
\snm{Sakai}\thanksref{t3}\corref{}\ead[label=e2]{sakai@cris.hokudai.ac.jp}
\ead[label=u1,url]{http://www.math.sci.hokudai.ac.jp/\textasciitilde sakai/}}
\runauthor{L.-C. Chen and A. Sakai}
\thankstext{t1}{Supported in part by the TJ\&MY Foundation and  NSC Grant 99-2115-M-030-004-MY3.}
\thankstext{t3}{Supported in part by the start-up fund of
L-Station at Hokkaido University and  JSPS Grant-in-Aid 21740059 for Young
Scientists (B).}

\affiliation{Fu-Jen Catholic University and Hokkaido University}

\address[A]{Department of Mathematics\\
Fu-Jen Catholic University\\
510 Chung Cheng Road\\ Hsinchuang\\
Taipei County 24205\\ Taiwan\\
\printead{e1}}

\address[B]{Creative Research Institution SOUSEI\\
Hokkaido University\\
North 21, West 10, Kita-ku\\
Sapporo 001-0021\\ Japan\\
\printead{e2}\\
}
\end{aug}

\received{\smonth{2} \syear{2010}}
\revised{\smonth{4} \syear{2010}}

%
\begin{abstract}
We consider random walk and self-avoiding walk whose 1-step
distribution is
given by $D$, and oriented percolation whose bond-occupation
probability is
proportional to $D$. Suppose that $D(x)$ decays as $|x|^{-d-\alpha}$ with
$\alpha>0$. For random walk in any dimension $d$ and for self-avoiding walk
and critical/subcritical oriented percolation above the common upper-critical
dimension $d_{\mathrm c}\equiv2(\alpha\wedge2)$, we prove large-$t$
asymptotics of the
gyration radius, which is the average end-to-end distance of random
walk/self-avoiding walk of length $t$ or the average spatial size of an
oriented percolation cluster at time~$t$. This proves the conjecture for
long-range self-avoiding walk in [\textit{Ann. Inst. H.
Poincar\'e Probab. Statist.} (2010), to appear] and for
long-range oriented
percolation in [\textit{Probab. Theory Related Fields} \textbf
{142} (2008) 151--188] and [\textit{Probab. Theory Related Fields}
\textbf{145} (2009) 435--458].
\end{abstract}



\end{frontmatter}

\section{Introduction}

\subsection{Motivation}

Let $\varphi_t^\mathrm{RW}(x)$ be the $t$-step transition
probability for
random walk
on ${\mathbb Z}^d$: $\varphi_0^{\mathrm{RW}}(x)=\delta_{o,x}$ and
%
\begin{eqnarray}\label{eq:RW}
\varphi_t^{\mathrm{RW}}(x)=(\varphi_{t-1}^
{\mathrm{RW}}*D)(x)\equiv\sum_{y\in{\mathbb Z}^d
}\varphi_{t-1}^{\mathrm{RW}}
(y) D(x-y)\qquad[t\in\mathbb{N}].
\end{eqnarray}
Suppose that the 1-step distribution $D$ is ${\mathbb Z}^d$-symmetric.
How does the
$r$th moment $\sum_x|x|^r\varphi_t^{\mathrm{RW}}(x)$ grow as
$t\to\infty$, where $|\cdot|$ denotes the Euclidean distance?
When $r=2$ and $\sigma^2\equiv\sum_x|x|^2D(x)<\infty$, the answer
is trivial:
$\sum_x|x|^2\varphi_t^{\mathrm{RW}}(x)=\sigma^2t$ since the
variance of the
sum of
independent random variables is the sum of their variances. It is not
so hard
to see that $\sum_x|x|^r\varphi_t^{\mathrm{RW}}(x)=O(t^{r/2})$
as $t\to\infty
$ for other
values of $r>2$, as long as $\sum_x|x|^rD(x)<\infty$. Even so, it
may not be that easy to identify the constant $C\in(0,\infty)$ such that
$(\sum_x|x|^r\varphi_t^{\mathrm{RW}}(x))^{1/r}\sim C\sqrt t$.
Here, and in the
rest of
the paper, ``$f(z)=O(g(z))$'' means that $|f(z)/g(z)|$ is bounded for
all $z$
in some relevant set, while ``$f(z)\sim g(z)$'' means that $f(z)/g(z)$ tends
to 1 in some relevant limit for $z$.

Let $\alpha>0$, $L\in[1,\infty)$ and suppose that
$D(x)\approx|x/L|^{-d-\alpha}$ for large $x$ such that its Fourier transform
$\hat D(k)\equiv\sum_{x\in{\mathbb Z}^d}e^{ik\cdot x}D(x)$ satisfies
%
\begin{eqnarray}\label{eq:D-asympt}
1-\hat D(k)=v_\alpha|k|^{\alpha\wedge2}\times
\cases{
1+O((L|k|)^\epsilon),&\quad$\alpha\ne2$,
\cr
\log\dfrac1{L|k|}+O(1),&\quad$\alpha=2$
}
\end{eqnarray}
for some $v_\alpha=O(L^{\alpha\wedge2})$ and $\epsilon>0$. If
$\alpha>2$ (or
$D$ is finite-range), then $v_\alpha\equiv\sigma^2/(2d)$. As shown in
Appendix~\hyperref[appendix:D]{A.1}, the long-range Kac potential
%
\begin{eqnarray}\label{eq:kac}
D(x)=\frac{h(y/L)}{\sum_{y\in{\mathbb Z}^d}h(y/L)}\qquad[x\in
{\mathbb Z}^d],
\end{eqnarray}
defined in terms of a rotation-invariant
function $h$ satisfying
\begin{eqnarray*}
h(x)=\frac{1+O((|x|\vee1)^{-\rho})}{(|x|\vee1)^{d+\alpha
}}\qquad
[x\in{\mathbb R}^d]
\end{eqnarray*}
for some $\rho>\epsilon$, satisfies the above properties. Notice that
$\sum_x|x|^rD(x)=\infty$ for $r\ge\alpha$ and, in particular,
$\sigma^2=\infty$ if $\alpha\le2$. This is of interest in investigating
the asymptotic behavior of $\sum_x|x|^r\varphi_t^\mathrm
{RW}(x)$ for \emph{all}
$r\in(0,\alpha)$ and understanding its $\alpha$-dependence.

In fact, our main interest is in proving sharp asymptotics of the gyration
radius of order $r\in(0,\alpha)$, defined as
\begin{eqnarray*}
\xi_t^{(r)}=\biggl(\frac{\sum_{x\in{\mathbb Z}^d}|x|^r\varphi
_t(x)}{\sum_{x\in{\mathbb Z}^d}
\varphi_t(x)}\biggr)^{1/r},
\end{eqnarray*}
where $\varphi_t(x)\equiv\varphi_t^\mathrm{SAW}(x)$ is the two-point
function for
$t$-step self-avoiding walk whose 1-step distribution is given by $D$, or
$\varphi_t(x)\equiv\varphi_t^\mathrm{OP}(x)$ is the two-point
function for oriented
percolation whose bond-occupation probability for each bond $((u,s),(v,s+1))$
is given by $pD(v-u)$, independently of $s\in{\mathbb Z}_+$, where
$p\ge0$ is the
percolation parameter. More precisely,
%
\begin{eqnarray}\label{eq:2pt-def}
\varphi_t(x)=
\cases{
\displaystyle\varphi_t^{\mathrm{RW}}(x)\equiv\mathop{\sum_{\omega:o\to x}}_{
(|\omega|=t)}
\prod_{s=1}^tD(\omega_s-\omega_{s-1}),
\cr
\displaystyle\varphi_t^\mathrm{SAW}(x)\equiv\mathop{\sum_{\omega:o\to x}}_
{(|\omega|=t)}
\prod_{s=1}^tD(\omega_s-\omega_{s-1})\prod_{0\le i<j\le t}(1-
\delta_{\omega_i,\omega_j}),
\cr
\varphi_t^\mathrm{OP}(x)\equiv{\mathbb P}_p((o,0)\to
(x,t)),
}
\end{eqnarray}
where $\prod_{0\le i<j\le t}(1-\delta_{\omega_i,\omega_j})$ is the
self-avoiding constraint on $\omega$ and $\{(o,0)\to(x,t)\}$ is the event
that either $(x,t)=(o,0)$ or there is a consecutive sequence of
occupied bonds
from $(o,0)$ to $(x,t)$ in the time-increasing direction. The gyration radius
$\xi_t^{(r)}$ represents a typical end-to-end
distance of a
linear structure of length $t$ or a typical spatial size of a cluster
at time
$t$. It has been expected (and would certainly be true for random walk
in any dimension) that, above the common upper-critical dimension
$d_\mathrm{c}=2(\alpha\wedge2)$ for self-avoiding walk and oriented
percolation,
for every $r\in(0,\alpha)$,
%
\begin{eqnarray}\label{eq:conjecture}
\xi_t^{(r)}=
\cases{
O\bigl(t^{1/({\alpha\wedge2})}\bigr),&\quad$\alpha\ne2$,
\cr
O\bigl(\sqrt{t\log t}\bigr),&\quad$\alpha=2$.
}
\end{eqnarray}
Heydenreich~\cite{hhh} proved (\ref{eq:conjecture}) for self-avoiding
walk, but
only for small $r<\alpha\wedge2$. Nevertheless, this small-$r$ result is
enough to prove weak convergence of self-avoiding walk to an $\alpha$-stable
process/Brownian motion, depending on the value of~$\alpha$ \cite{hhh}.

As stated below in Theorem~\ref{theorem:main2}, we prove sharp asymptotics
(including the proportionality constant) of
$\sum_x|x_1|^r\varphi_t(x)/\sum_x\varphi_t(x)$ as $t\to\infty$,
where $x_1$
is the first coordinate of $x\equiv(x_1,\dots,x_d)$, and show that
(\ref{eq:conjecture}) holds for all $r\in(0,\alpha)$, solving the
open problems
in \cite{csII,hhh}.

\subsection{Main results}\label{ss:main}
Let $m_{\rm c}\ge1$ be the model-dependent radius of convergence for
the sequence
$\sum_x\varphi_t(x)$. For random walk, $m_{\rm c}=1$ since $\sum
_x\varphi_t^{\mathrm{RW}}(x)$
is always 1. For self-avoiding walk, $m_{\rm c}>1$ due to the
self-avoiding constraint in (\ref{eq:2pt-def}) and, indeed, $m_{\rm
c}=1+O(L^{-d})$ for
$d>d_{\mathrm c}$ and $L\gg1$ \cite{hhs08}. For oriented percolation,
$m_{\rm
c}$ depends on
the percolation parameter $p$ [i.e., $m_{\rm c}=m_{\rm c}(p)$] and was
denoted by $m_p$ in
\cite{csI,csII}. It has been proven \cite{csI} that $m_{\rm c}(p)>1$
for $p<p_{\mathrm c},$
and $m_{\rm c}(p_{\mathrm c})=1$ for $d>d_{\mathrm c}$ and $L\gg1$,
where $p_{\mathrm c}$ is the
critical point
characterized by the divergence of the susceptibility:
$\sum_{t=0}^\infty\sum_{x\in{\mathbb Z}^d}\varphi_t^\mathrm
{OP}(x)\uparrow\infty
$ as
$p\uparrow p_{\mathrm c}$. It has also been\vspace*{1pt} proven \cite{csI} that
$pm_\mathrm{c}=1+O(L^{-d})$
for all $p\le p_{\mathrm c}$.\vspace*{1pt}

Let $C_{\mathrm I}$ and $C_{\mathrm{II}}$ be the constants in \cite
{csI,csII,hhh} such that,
as $t\to\infty$,
%
\begin{eqnarray}\label{eq:main-I&II}
\sum_{x\in{\mathbb Z}^d}\varphi_t(x)\sim C_{\mathrm I} m_{\rm c}^{-t},
     \qquad
\frac{\sum_{x\in{\mathbb Z}^d}e^{ik_t\cdot x}\varphi_t(x)}{\sum
_{x\in{\mathbb Z}^d
}\varphi_t(x)}
\sim e^{-C_{\mathrm{II}}|k|^{\alpha\wedge2}},
\end{eqnarray}
where
%
\begin{eqnarray}\label{eq:k-scaling}
k_t=k\times
\cases{
(v_\alpha t)^{-1/({\alpha\wedge2})},&\quad$\alpha\ne2$,
\cr
\bigl(v_2 t\log\sqrt t\bigr)^{-1/2},&\quad$\alpha=2$.
}
\end{eqnarray}
Because of this scaling, we have $C_{\mathrm I}^\mathrm
{RW}=C_{\mathrm{II}}^{\mathrm{RW}}=1$ for random
walk. For
self-avoiding walk and critical/subcritical oriented percolation for
$d>2(\alpha\wedge2)$ with $L\gg1$ (depending on the models), it has been
proven that the model-dependent constants $C_{\mathrm I}$ and
$C_{\mathrm{II}}$ are both
$1+O(L^{-d})$ \cite{csI,hhh} and that the $O(L^{-d})$ term in
$C_{\mathrm{II}}$ exhibits
crossover behavior at $\alpha=2$ \cite{csII,hhh}. We will provide precise
expressions for $C_{\mathrm I}$ and $C_{\mathrm{II}}$ at the end of
Section~\ref{ss:outline}.

Our first result is the following asymptotic behavior of the generating
function for the sequence $\sum_x|x_1|^r\varphi_t(x)$.

\begin{theorem}\label{theorem:main1}
Consider the three aforementioned long-range models. For random walk in any
dimension $d$ with any $L$, and for self-avoiding walk and critical/subcritical
oriented percolation for $d>d_{\mathrm c}\equiv2(\alpha\wedge2)$
with $L\gg
1$ (depending
on the models), the following holds for all $r\in(0,\alpha)$: as
$m\uparrow m_{\rm c}$,
\begin{eqnarray}\label{eq:main1}
\qquad\sum_{t=0}^\infty m^t\sum_{x\in{\mathbb Z}^d}|x_1|^r\varphi
_t(x)&=&\frac
{2\sin
({r\pi}/{(\alpha\vee2))}}{(\alpha\wedge2)\sin({r\pi}/\alpha)} \Gamma(r+1)
 \frac{C_{\mathrm I}(C_{\mathrm{II}}v_\alpha)^{{r}/{(\alpha
\wedge2)}}}{(1-
{m}/{m_{\rm c}})^{1+
{r}/{(\alpha\wedge2)}}}\nonumber
\\[-8pt]\\[-8pt]
&&{}\times
\cases{
1+O\biggl(\biggl(1-\dfrac{m}{m_{\rm c}}\biggr)^\epsilon\biggr),&\quad$\alpha\ne2$,
\cr
\biggl(\log\dfrac1{\sqrt{1-{m}/{m_{\rm c}}}}
\biggr)^{r/2}+O(1),&\quad$\alpha=2$
}
\nonumber
\end{eqnarray}
for some $\epsilon>0$ when $\alpha\ne2$. The $O(1)$ term for $\alpha=2$
is independent of $m$.
\end{theorem}

It is worth emphasizing that, although $C_{\mathrm I},C_{\mathrm
{II}},m_{\rm c}$ are
model-dependent,
the formula (\ref{eq:main1}) itself is universal. Expanding (\ref
{eq:main1}) in
powers of $m$ and using (\ref{eq:main-I&II}), we obtain the following theorem.

\begin{theorem}\label{theorem:main2}
Under the same condition as in Theorem~\ref{theorem:main1}, as $t\to
\infty$,
%
\begin{eqnarray}\label{eq:main2}
\frac{\sum_{x\in{\mathbb Z}^d}|x_1|^r\varphi_t(x)}{\sum_{x\in
{\mathbb Z}^d}\varphi
_t(x)}&\sim&
\frac{2\sin({r\pi/(\alpha\vee2)})}{(\alpha\wedge2)\sin
({r\pi}/\alpha)}
\frac{\Gamma(r+1)}{\Gamma({r}/{(\alpha\wedge2)}+1)}\nonumber
\\[-8pt]\\[-8pt]
&&{}\times
\cases{
\bigl(C_{\mathrm{II}}v_\alpha t\bigr)^{{r}/{(\alpha\wedge2)}},&\quad$\alpha\ne2$,
\cr
\bigl(C_{\mathrm{II}}v_2t\log\sqrt{t}\bigr)^{r/2},&\quad$\alpha=2$.
}
\nonumber
\end{eqnarray}
\end{theorem}

We note that $C_{\mathrm{II}}$ is the only model-dependent term in (\ref
{eq:main2}).
As far as we are aware, the sharp asymptotics (\ref{eq:main1}) and (\ref
{eq:main2})
for \emph{all} real $r\in(0,\alpha)$ are new, even for random walk.

Although we focus our attention on the long-range models defined by $D$ that
satisfies (\ref{eq:D-asympt}), our proof also applies to
finite-range models, for which $\alpha$ is considered to be infinity.

Using $|x_1|^r\le|x|^r\le d^{r/2}\sum_{j=1}^d|x_j|^r$ and the
${\mathbb Z}^d$-symmetry
of the models, we are finally able to arrive at the following result.

\begin{corollary}
Under the same condition as in Theorem~\ref{theorem:main1}, (\ref
{eq:conjecture})
holds for all $r\in(0,\alpha)$. In particular, when $r=2<\alpha$,
%
\begin{eqnarray}\label{eq:mean-square}
\xi_t^{(2)}\mathop{\sim}_{t\to\infty}\sqrt{C_{\mathrm{II}}\sigma^2 t}.
\end{eqnarray}
\end{corollary}

As mentioned earlier, (\ref{eq:conjecture}) has been proven \cite
{hhh} for
self-avoiding walk, but only for small $r<\alpha\wedge2$. The sharp
asymptotics (\ref{eq:mean-square}) has been proven \cite{hs02} for
self-avoiding
walk and critical oriented percolation defined by $D$ that has a finite
$(2+\epsilon)$th moment for some $\epsilon>0$. Our proof is based on a
different method than those used in \cite{hhh,hs02}. It is closer to the
method, explained in the next subsection, used in \cite{ms93} for finite-range
self-avoiding walk and in \cite{ny95} for critical/subcritical finite-range
oriented percolation.

We strongly believe that the same method should work
for lattice trees. Any two points in a lattice tree are connected by a
unique path, so the number of bonds contained in that path can be considered
as time and we can apply the current method to obtain the same
results (with different values for $C_{\mathrm I},C_{\mathrm{II}}$). As
this suggests,
time, or
something equivalent, is important for the current method to work. For
unoriented percolation, for example, it is not so clear what should be
interpreted as time. However, if $D$ is biased in average in one direction,
say, the positive direction of the first coordinate, then $x_1$ can be treated
as time and, after subtracting the effect of the bias, we may obtain the
results even for unoriented percolation.

\subsection{Outline and notation}\label{ss:outline}
In this subsection, we outline the proof of Theorem~\ref
{theorem:main1} and
introduce some notation which is used in the rest of the paper. We also refer
interested readers to an extended version of this subsection in \cite{s09}.

One of the key elements for the proof is to represent the left-hand
side of
(\ref{eq:main1}) in terms of the generating function (i.e., the
Fourier--Laplace
transform)
of the two-point function. We now explain this representation.

Given a function $f_t(x)$, where
$(x,t)\in{\mathbb Z}^d\times{\mathbb Z}_+$, we formally define
\begin{eqnarray*}
\hat f(k,m)=\sum_{t=0}^\infty m^t\sum_{x\in{\mathbb Z}^d}f_t(x)
e^{ik\cdot
x}\qquad[k\in
[-\pi,\pi]^d,~m\ge0].
\end{eqnarray*}
We note that $\hat\varphi(k,m)$ is well defined when $m<m_{\rm c}$
(recall that
$m_{\rm c}\ge1$, as explained at the beginning of Section~\ref{ss:main}).
Let\vspace*{-2pt}
%
\begin{eqnarray}\label{eq:nablone}
\nabla_{1}^n\hat f(l,m)=\frac{\partial^n\hat f(k,m)}{\partial
k_1^n}\bigg|_{k=l}
\qquad[l\in[-\pi,\pi]^d,~n\in{\mathbb Z}_+].\vspace*{-2pt}
\end{eqnarray}
Then, for $r=2j<\alpha$ ($j\in\mathbb{N}$), we obtain the representation\vspace*{-2pt}
\begin{eqnarray*}
\sum_{t=0}^\infty m^t\sum_{x\in{\mathbb
Z}^d}x_1^{2j}f_t(x)=(-1)^j\nabla_{1}
^{2j}\hat
f(0,m).\vspace*{-2pt}
\end{eqnarray*}

For $r\in(0,\alpha\wedge2)$, we generate the factor $|x_1|^r$ by
using the
constant $K_r\in(0,\infty)$, as follows (see \cite{csII}):\vspace*{-2pt}
%
\begin{eqnarray}\label{eq:Kr}
K_r\equiv\int_0^\infty\frac{1-\cos v}{v^{1+r}} \,\mathrm{d}v=|x_1|^{-r}
\int_0^\infty\frac{1-\cos(ux_1)}{u^{1+r}} \,\mathrm{d}u.\vspace*{-2pt}
\end{eqnarray}
Suppose, from now on, that $f_t$ is ${\mathbb Z}^d$-symmetric. Then,\vspace*{-2pt}
\begin{eqnarray*}
\sum_{t=0}^\infty m^t\sum_{x\in{\mathbb Z}^d}|x_1|^rf_t(x)&=&\frac
1{K_r}\int
_0^\infty
\frac{\mathrm{d}u}{u^{1+r}}\sum_{t=0}^\infty m^t\sum_{x\in{\mathbb
Z}^d}
\bigl(1-\cos(ux_1)
\bigr)f_t(x)
\\[-2pt]
&=&\frac1{K_r}\int_0^\infty\frac{\mathrm{d}u}{u^{1+r}} \bigl(\hat
f(0,m)-\hat
f(\vec u,m)\bigr),\vadjust{\goodbreak}
\end{eqnarray*}

\noindent
where $\vec u=(u,0,\dots,0)\in{\mathbb R}^d$. Let
\begin{eqnarray}\label{eq:laplace}
\bar\Delta_l\hat f(k,m)
&\equiv&\hat f(k,m)-\frac{\hat f(k+l,m)+\hat f(k-l,m)}2\nonumber
\\[-8pt]\\[-8pt]
&=&\sum_{t=0}^\infty m^t\sum_{x\in{\mathbb Z}^d}\bigl(1-\cos(l\cdot
x)
\bigr)f_t(x)
e^{ik\cdot x}.\nonumber
\end{eqnarray}
We note that $\bar\Delta_l\hat f(k,m)$ is equivalent to
$\frac{-1}2\Delta_l\hat f(k,m)$ in the previous papers (e.g.,
\cite{csI,csII}). In particular,
\begin{eqnarray*}
\bar\Delta_l\hat f(0,m)=\hat f(0,m)-\hat f(l,m).
\end{eqnarray*}
Therefore, for $r\in(0,\alpha\wedge2)$,
\begin{eqnarray*}
\sum_{t=0}^\infty m^t\sum_{x\in\mathbb{Z}^d}|x_1|^rf_t(x)=\frac1{K_r}
\int_0^\infty\frac{\mathrm{d}u}{u^{1+r}} \bar\Delta_{\vec u}\hat f(0,m).
\end{eqnarray*}
For $r=2j+q<\alpha$ [$j\in\mathbb{N}$, $q\in(0,2)$], we combine the above
representations as
\begin{eqnarray*}
\sum_{t=0}^\infty m^t\sum_{x\in{\mathbb
Z}^d}|x_1|^{2j+q}f_t(x)&=&\frac1{K_q}
\int_0^\infty\frac{\mathrm{d}u}{u^{1+q}}\sum_{t=0}^\infty m^t\sum
_{x\in{\mathbb Z}^d}
\bigl(1-\cos(ux_1)\bigr)x_1^{2j}f_t(x)
\\
&=&\frac{(-1)^j}{K_q}\int_0^\infty\frac{\mathrm{d}u}{u^{1+q}}
\bigl(\nabla_{1}^{2j}
\hat f(0,m)-\nabla_{1}^{2j}\hat f(\vec u,m)\bigr)
\\
&=&\frac{(-1)^j}{K_q}\int_0^\infty\frac{\mathrm{d}u}{u^{1+q}}
\bar\Delta_{\vec u}
\nabla_{1}^{2j}\hat f(0,m).
\end{eqnarray*}

From now on, as long as no confusion arises, we will simply omit $m$ and
abbreviate $\hat f(k,m)$ to $\hat f(k)$. Then, the above three
representations are summarized as
\begin{eqnarray}\label{eq:representation}
&&\sum_{t=0}^\infty m^t\sum_{x\in{\mathbb Z}^d}|x_1|^rf_t(x)\nonumber
\\[-8pt]\\[-8pt]
&&\qquad =
\cases{
(-1)^j\nabla_{1}^{2j}\hat f(0)\qquad [r=2j<\alpha,~j\in\mathbb{N}],
\cr
\displaystyle\dfrac{(-1)^j}{K_q}\int_0^\infty\dfrac{\mathrm
{d}u}{u^{1+q}}
\bar\Delta_{\vec u}\nabla_{1}^{2j}\hat f(0)
\cr
\qquad [r=2j+q<\alpha, j
\in{\mathbb Z}_+, q\in(0,2)].
}\nonumber
\end{eqnarray}
Also, we will abbreviate $\hat f(k,m_{\rm c})$ to $\hat f_{\mathrm c}(k)$
whenever it
is well defined. Moreover, we will use the notation
\begin{eqnarray*}
\partial_m\hat f_{\mathrm c}(k)=\frac{\partial\hat f(k,m)}{\partial
m}\bigg|_{m
\uparrow m_{\rm c}}.
\end{eqnarray*}

Another key element for the proof of the main theorem is the lace expansion
(see, e.g., \cite{s06}, Sections~3 and 13),
%
\begin{eqnarray}\label{eq:lace-exp}
\varphi_t(x)=I_t(x)+\sum_{s=1}^t(J_s*\varphi_{t-s})(x),
\end{eqnarray}
where, for $t\ge0$,
%
\begin{eqnarray}\label{eq:I-def}
I_t(x)=
\cases{
\delta_{x,o}\delta_{t,0},&\quad\mbox{RW/SAW},
\cr
\pi_t^\mathrm{OP}(x),&\quad$\mathrm{OP}$,
}
\end{eqnarray}
and for $t\ge1$,
%
\begin{eqnarray}\label{eq:J-def}
J_t(x)=
\cases{
D(x)\delta_{t,1},&\quad\mbox{RW},
\cr
D(x)\delta_{t,1}+\pi_t^\mathrm{SAW}(x),&\quad\mbox{SAW},
\cr
(\pi_{t-1}^\mathrm{OP}*pD)(x),&\quad\mbox{OP}.
}
\end{eqnarray}
Recall (\ref{eq:RW}) for random walk. For self-avoiding walk and oriented
percolation, $\pi_t^\mathrm{SAW}(x)$ and $\pi_t^\mathrm
{OP}(x)$ are (alternating
sums of)
the model-dependent lace expansion coefficients (see, e.g., \cite
{s06} for
their precise definitions). By~(\ref{eq:lace-exp}), we obtain
%
\begin{eqnarray}\label{eq:lace-exp-fourier}
\hat\varphi(k)=\hat I(k)+\hat J(k) \hat\varphi(k).
\end{eqnarray}
From this, we can derive identities for the ``derivatives'' of
$\hat\varphi$ in (\ref{eq:representation}). For example,
%
\begin{eqnarray}\label{eq:derivative-example}
\bar\Delta_{\vec u}\hat\varphi(0)\equiv\hat\varphi(0)-\hat
\varphi
(\vec u)
&=&\hat I(0)+\hat J(0) \hat\varphi(0)-\bigl(\hat I(\vec u)+\hat
J(\vec u)
\hat\varphi(\vec u)\bigr)\nonumber
\\
&=&\bar\Delta_{\vec u}\hat I(0)+\hat J(0) \hat\varphi(0)-\hat
J(\vec
u) \hat
\varphi(\vec u)\nonumber
\\[-8pt]\\[-8pt]
&=&\bar\Delta_{\vec u}\hat I(0)+\hat\varphi(0) \bar\Delta_{\vec
u}\hat
J(0)+\hat
J(\vec u) \bar\Delta_{\vec u}\hat\varphi(0)\nonumber
\\
&=&\frac1{1-\hat J(\vec u)}\bigl(\bar\Delta_{\vec u}\hat I(0)+\hat
\varphi(0)
\bar\Delta_{\vec u}\hat J(0)\bigr),\nonumber
\end{eqnarray}
where the last line has been obtained by solving the previous equation for
$\bar\Delta_{\vec u}\hat\varphi(0)$. Hence, for $r\in(0,\alpha
\wedge2)$,
\begin{eqnarray}\label{eq:derivative-appl}
\sum_{t=0}^\infty m^t\sum_{x\in{\mathbb Z}^d}|x_1|^r\varphi
_t(x)&=&\frac
{\hat\varphi(0)}
{K_r}\int_0^\infty\frac{\mathrm{d}u}{u^{1+r}} \frac{\bar\Delta
_{\vec
u}\hat J
(0)}{1-\hat J(\vec u)}\nonumber
\\[-8pt]\\[-8pt]
&&{}+\frac1{K_r}\int_0^\infty\frac{\mathrm{d}u}{u^{1+r}}
\frac{\bar\Delta_{\vec u}\hat I(0)}{1-\hat J(\vec u)}.\nonumber
\end{eqnarray}
It is known \cite{csI,hhs08} that as long as $d>d_{\mathrm c}$ (and
$L\gg1$), it
is easier to tame $\hat I$ and $\hat J$, up to $m=m_{\rm c}$, than to tame
$\hat\varphi$. We will thus be able to analyze the integrals on the
right-hand side of (\ref{eq:derivative-appl}) and prove the main theorem.

Before closing this subsection, we provide the following
representations for
the constants $C_{\mathrm I}$ and $C_{\mathrm{II}}$ in (\ref{eq:main1})
in terms of $\hat
I_{\mathrm c}$ and
$\hat J_{\mathrm c}$:
%
\begin{eqnarray}\label{eq:CICII}
C_{\mathrm I}=\frac{\hat I_{\mathrm c}(0)}{m_{\rm c}\,\partial_m\hat
J_{\mathrm c}(0)},\qquad
C_{\mathrm{II}}=\frac1{m_{\rm c}\,\partial_m\hat J_{\mathrm c}(0)}\lim
_{k\to0}\frac
{\bar\Delta_k\hat J_{\mathrm c}
(0)}{\bar\Delta_k\hat D(0)}.
\end{eqnarray}
In Section~\ref{s:CICII}, we will explain the heuristics for the derivation
of these representations.

\subsection{Organization}
In the remainder of the paper, whenever we consider self-avoiding walk
and oriented
percolation, we assume $d>d_{\mathrm c}$ and $L\gg1$, as well as $p\le
p_{\mathrm c}$
for oriented
percolation.

The paper is organized as follows. In Section~\ref{s:proof}, we
prove Theorem~\ref{theorem:main1} for $r\in(0,\alpha\wedge2)$
(Section~\ref{ss:0<r<alpha,2}), for $r=2j<\alpha$ with $j\in\mathbb{N}$
(Section~\ref{ss:r=2j}) and for $r=2j+q<\alpha$ with $j\in\mathbb
{N}$ and
$q\in(0,2)$
(Section~\ref{ss:r=2j+q}) separately, assuming
Propositions~\ref{proposition:r<alpha&2} and \ref
{proposition:2<r<alpha}. We
prove those key propositions in Section~\ref{s:key-prop}.

We strongly believe that the results for self-avoiding walk and oriented
percolation are the most important and interesting parts of this work.
However, for those who are more interested in random walk,
we make the following suggestion: read up to Section~\ref{s:proof}
for the proof of Theorem~\ref{theorem:main1}, where
Propositions~\ref{proposition:r<alpha&2} and \ref{proposition:2<r<alpha}
are used. However, Proposition~\ref{proposition:r<alpha&2} and a part
[i.e., (\ref{eq:2<r<alpha-IJbds})] of Proposition~\ref{proposition:2<r<alpha}
are trivial for random walk. The remaining part [i.e.,
(\ref{eq:2<r<alpha-varphibd})] of Proposition~\ref{proposition:2<r<alpha}
is the result of Lemma~\ref{lmm:induction1}, which is proved in
Section~\ref{ss:prop2}.

\section{Preliminaries}\label{s:CICII}
In this section, we review in outline the derivation in \cite
{csI,csII,hhh} of
the constants $C_{\mathrm I}$ and $C_{\mathrm{II}}$. During the course
of this, we
summarize the
already known properties of $\hat I$ and $\hat J$ and introduce some
quantities used in the following sections.

We begin by solving (\ref{eq:lace-exp-fourier}) for $\hat\varphi(k)$,
which yields
%
\begin{eqnarray}\label{eq:lace-exp-solution}
\hat\varphi(k)=\frac{\hat I(k)}{1-\hat J(k)},
\end{eqnarray}
where, by (\ref{eq:I-def}) and (\ref{eq:J-def}),
%
\begin{eqnarray}\label{eq:Ihat}
\hat I(k)&=&
\cases{
1,&\quad\mbox{RW/SAW},
\cr
\hat\pi^\mathrm{OP}(k),&\quad\mbox{OP},
}
\\\label{eq:Jhat}
\hat J(k)&=&
\cases{
m\hat D(k),&\quad\mbox{RW},
\cr
m\hat D(k)+\hat\pi^\mathrm{SAW}(k),&\quad\mbox{SAW},
\cr
\hat\pi^\mathrm{OP}(k)pm\hat D(k),&\quad\mbox{OP}.
}
\end{eqnarray}
It is known \cite{csI,hhs08} that
%
\begin{eqnarray}\label{eq:property1}
\hat\pi^\mathrm{SAW}(k)=O(L^{-d}) m^2,\qquad
\hat\pi^\mathrm{OP}(k)-1=O(L^{-d}) (pm)^2,
\end{eqnarray}
where the $O(L^{-d})$ terms are uniform in $k\in[-\pi,\pi]^d$ and
$m\le m_{\rm c}$.
Therefore, $\hat I(k)$ and $\hat J(k)$ are both convergent for all
$k\in[-\pi,\pi]^d$ and $m\le m_{\rm c}$. However, since
$\hat\varphi(0)$ diverges as $m\uparrow m_{\rm c}$, we can
characterize $m_{\rm c}$
by the equation
%
\begin{eqnarray}\label{eq:critpt}
1=\hat J_{\mathrm c}(0)=
\cases{
m_{\rm c},&\quad\mbox{RW},
\cr
m_{\rm c}+\hat\pi_{\mathrm c}^\mathrm{SAW}(0),&\quad\mbox
{SAW},
\cr
\hat\pi_{\mathrm c}^\mathrm{OP}(0)pm_{\rm c},&\quad\mbox{OP}.
}
\end{eqnarray}
Using this identity, we obtain that, as $m\uparrow m_{\rm c}$ (see
\cite{csI,hhh} for
the precise argument),
\begin{eqnarray*}
\hat\varphi(k)&=&\frac{\hat I(k)}{\hat J_{\mathrm c}(0)-\hat
J_{\mathrm c}
(k)+m_{\rm c}{((\hat
J_{\mathrm c}(k)-\hat J(k))/(m_{\rm c}-m))} (1-{m}/{m_{\rm c}})}
\\
&\sim&\frac{\hat I_{\mathrm c}(k)}{\bar\Delta_k\hat J_{\mathrm
c}(0)+m_{\rm
c}\,\partial_m\hat J_{\mathrm c}(k)
(1-{m/m_{\rm c}})}
\\
&=&\frac{\hat I_{\mathrm c}(k)}{\bar\Delta_k\hat J_{\mathrm
c}(0)+m_{\rm c}\,\partial
_m\hat J_{\mathrm c}(k)}
\sum_{t=0}^\infty\biggl(\frac{m_{\rm c}\,\partial_m\hat J_{\mathrm c}
(k)}{\bar\Delta_k\hat
J_{\mathrm c}(0)+m_{\rm c}\,\partial_m\hat J_{\mathrm c}(k)} \frac
{m}{m_{\rm
c}}\biggr)^t,
\end{eqnarray*}
hence,
%
\begin{eqnarray}\label{eq:lonkey-markus}
\sum_{x\in{\mathbb Z}^d}\varphi_t(x) e^{ik\cdot x}&\mathop{\sim}\limits_{t\to
\infty}&
\frac{\hat I_{\mathrm c}(k)}{\bar\Delta_k\hat J_{\mathrm
c}(0)+m_{\rm c}\,\partial
_m\hat
J_{\mathrm c}(k)} m_{\rm c}^{-t}\nonumber
\\[-8pt]\\[-8pt]
&&{}\times\biggl(1-\frac{\bar\Delta_k\hat J_{\mathrm
c}(0)}{\bar\Delta
_k\hat J_{\mathrm c}(0)
+m_{\rm c}\,\partial_m\hat J_{\mathrm c}(k)}\biggr)^t.\nonumber
\end{eqnarray}
In particular, at $k=0$,
%
\begin{eqnarray}\label{eq:lonkey-markus0}
\sum_{x\in{\mathbb Z}^d}\varphi_t(x)\sim\frac{\hat I_{\mathrm
c}(0)}{m_{\rm
c}\,\partial_m
\hat J_{\mathrm c}(0)} m_{\rm c}^{-t},
\end{eqnarray}
which yields the representation for $C_{\mathrm I}$ in (\ref{eq:CICII}).

In the above computation, we have used the fact that the quantities
such\vspace*{1pt} as
$m_{\rm c}\,\partial_m\hat J_{\mathrm c}(0)$ and $\bar\Delta_k\hat
J_{\mathrm c}(0)$
are all convergent
uniformly in $k\in[-\pi,\pi]^d$. To see this, we note that, by
(\ref{eq:Jhat}),
%
\begin{eqnarray}
\qquad m_{\rm c}\partial_m\hat J_{\mathrm c}(k)&=&
\cases{
m_{\rm c}\hat D(k),&\quad\mbox{RW},
\cr
m_{\rm c}\hat D(k)+m_{\rm c}\,\partial_m\hat\pi_{\mathrm
c}^\mathrm{SAW}(k),&\quad\mbox
{SAW},
\cr
\bigl(\hat\pi_{\mathrm c}^\mathrm{OP}(k)+m_{\rm c}\,\partial
_m\hat\pi_{\mathrm c}^\mathrm{OP}
(k)\bigr)pm_{\rm c}\hat D(k),
&\quad\mbox{OP},
}
\\\label{eq:Jhat-laplace}
\bar\Delta_k\hat J(0)&=&
\cases{
m\bar\Delta_k\hat D(0),&\quad\mbox{RW},
\cr
m\bar\Delta_k\hat D(0)+\bar\Delta_k\hat\pi^\mathrm
{SAW}(0),&\quad\mbox{SAW},
\cr
\bigl(\hat\pi^\mathrm{OP}(0) \bar\Delta_k\hat D(0)+\hat
D(k) \bar\Delta
_k\hat\pi^\mathrm{OP}(0)
\bigr)pm,&\quad\mbox{OP}.
}
\end{eqnarray}
However, it is known that $\pi^\mathrm{SAW}$ and $\pi
^\mathrm{OP}$ both satisfy
%
\begin{eqnarray}\label{eq:property2}
\qquad \quad|m_{\rm c}\,\partial_m\hat\pi_{\mathrm c}(k)|&\le&\sum_{t=0}^\infty t
m_{\rm
c}^t\sum_{x\in{\mathbb Z}^d}|\pi_t
(x)|\le O(L^{-d}),
\\\label{eq:property3}
|\bar\Delta_k\hat\pi(0)|&\le&\sum_{t=0}^\infty m_{\rm c}^t\sum
_{x\in
{\mathbb Z}^d}\bigl(1-\cos(k
\cdot x)\bigr) |\pi_t(x)|\le O(L^{-d}) \bar\Delta_k\hat
D(0)
\end{eqnarray}
for all $k\in[-\pi,\pi]^d$ and $m\le m_{\rm c}$ for the latter (see
\cite{csII}, Proposition~1, \cite{hhh}, the paragraph below
Theorem~1.2 and
\cite{hhs08}, Proposition~4.1, with an improvement due to monotone
convergence). By these bounds and using (\ref{eq:property1}) and (\ref
{eq:critpt}) and
the fact that $m_{\rm c}^\mathrm{SAW}$ and $pm_{\rm
c}^\mathrm{OP}$ are both
$1+O(L^{-d})$ (see the
beginning of Section~\ref{ss:main}), we conclude that
$m_{\rm c}\,\partial_m\hat J_{\mathrm c}(0)=1+O(L^{-d})$ and
$\bar\Delta_k\hat J_{\mathrm c}(0)=O(\bar\Delta_k\hat D(0))$.

Moreover, it has been proven \cite{csI,csII,hhh} that there exist
$\epsilon=\epsilon(d,\alpha)>0$ and $\delta=\delta(d,\alpha)$,
which is zero
if $\alpha=2$ and $>0$ if $\alpha\ne2$, such that $\pi
^\mathrm{SAW}$ and
$\pi^\mathrm{OP}$ both
satisfy
\begin{eqnarray*}
\sum_{t=0}^\infty t^{1+\epsilon}m_{\rm c}^t\sum_{x\in{\mathbb
Z}^d}|\pi
_t(x)|<\infty,\qquad
\sum_{t=0}^\infty m_{\rm c}^t\sum_{x\in{\mathbb Z}^d}|x|^{\alpha
\wedge
2+\delta}|\pi_t(x)|
<\infty.
\end{eqnarray*}
These bounds imply (see \cite{csI}, equations (6.13) and (6.14),
\cite{csII}, equations (3.3)--(3.4), \cite{hhh}, equations (2.25)--(2.28) and (2.64)--(2.70))
%
\begin{eqnarray}\label{eq:Jhat-diff}
\frac{\hat J_{\mathrm c}(0)-\hat J(0)}{1-{m}/{m_{\rm c}}}&=&m_{\rm
c}\,\partial_m\hat J_{\mathrm c}(0)
+O\biggl(\biggl(1-\dfrac{m}{m_{\rm c}}\biggr)^\epsilon\biggr),
\\\label{eq:M-source}
\frac{\bar\Delta_k\hat J(0)}{\bar\Delta_k\hat D(0)}&=&M+
\cases{
O(|k|^\delta),&\quad$\alpha\ne2$,
\cr
O\biggl(1/\log\dfrac1{|k|}\biggr),&\quad$\alpha=2$,
}
\end{eqnarray}
where the error terms in (\ref{eq:M-source}), which are zero for
random walk,
are uniform in $m\le m_{\rm c}$ and where $M\equiv M(m)$ is defined as
%
\begin{eqnarray}\label{eq:M-def}
M=
\cases{
m,&\quad\mbox{RW},
\cr
\displaystyle m+\dfrac{\nabla_{1}^2\hat\pi^\mathrm
{SAW}(0)}{-2v_\alpha}{\mathbh1}_{\{\alpha>2\}},
&\quad\mbox{SAW},
\cr
\displaystyle\biggl(\hat\pi^\mathrm{OP}(0)+\dfrac{\nabla
_{1}^2\hat\pi^\mathrm{OP}
(0)}{-2v_\alpha}
{\mathbh1}_{\{\alpha>2\}}\biggr)pm,&\quad\mbox{OP}.
}
\end{eqnarray}
The crossover terms, which are proportional to ${\mathbh1}_{\{\alpha
>2\}}$,
converge for all $m\le m_{\rm c}$ \cite{csII,hhh}.
By~(\ref{eq:lonkey-markus}) and (\ref{eq:lonkey-markus0}) and (\ref
{eq:M-source}), and
using $\lim_{t\to\infty}t \bar\Delta_{k_t}\hat D(0)=|k|^{\alpha
\wedge2}$,
due to the scaling (\ref{eq:k-scaling}), we obtain that, as $t\to
\infty$,
\begin{eqnarray*}
&&\frac{\sum_{x\in{\mathbb Z}^d}\varphi_t(x) e^{ik_t\cdot x}}{\sum
_{x\in{\mathbb Z}^d
}\varphi_t(x)}
\\
&&\qquad \sim\biggl(1-\frac{\bar\Delta_{k_t}\hat J_{\mathrm c}(0)}{\bar
\Delta
_{k_t}\hat J_{\mathrm c}(0)+m_{\rm c}
\,\partial_m\hat J_{\mathrm c}(k_t)}\biggr)^t\qquad[\because(\ref{eq:lonkey-markus})_{k=k_t}
/(\ref{eq:lonkey-markus0})]
\\
&&\qquad \sim\exp\biggl(-\frac{\bar\Delta_{k_t}\hat J_{\mathrm
c}(0)}{m_{\rm
c}\,\partial_m\hat J_{\mathrm c}(0)}
 t\biggr)\qquad[\because 1-\tau\sim e^{-\tau}\mbox{ as }\tau\to
0]
\\
&&\qquad=\exp\biggl(-\frac1{m_{\rm c}\,\partial_m\hat J_{\mathrm c}(0)}
\frac
{\bar\Delta_{k_t}\hat J_{\mathrm c}
(0)}{\bar\Delta_{k_t}\hat D(0)} t \bar\Delta_{k_t}\hat D(0)
\biggr)
\\
&&\qquad\sim\exp\biggl(-\frac{M_{\mathrm c}}{m_{\rm c}\,\partial_m\hat
J_{\mathrm c}(0)}
|k|^{\alpha\wedge2}
\biggr),
\end{eqnarray*}
where $M_{\mathrm c}=M(m_{\rm c})$. This yields the representation for
$C_{\mathrm{II}}$ in
(\ref{eq:CICII}).

\begin{Remark*}
It is natural for some readers to wonder why we do not directly prove
(\ref{eq:main2}) by using the formula (\ref{eq:lonkey-markus}) for
$\sum_x\varphi_t(x) e^{ik\cdot x}$, instead of proving the asymptotics
(\ref{eq:main1}) of its generating function and expanding it in powers
of $m$.
In fact, the first-named author was able to derive an asymptotic
expression for $\sum_x|x_1|^r\varphi_t(x)$ using (\ref
{eq:lonkey-markus}), but
the proportionality constant was in a rather complicated sum form. We then
concluded that using (\ref{eq:lonkey-markus}) would not be an ideal
method for
deriving the simplest possible display of the proportionality constant and
started searching for another method. That turns out to be the use of
the generating
function, as explained in this paper. Later, the first-named author
proved that the aforementioned
sum form is indeed an expansion of the proportionality constant in
(\ref{eq:main2}).
\end{Remark*}

\section{Proof of the main results}\label{s:proof}

\subsection{\texorpdfstring{Proof of Theorem~\protect\ref{theorem:main1} for $r\in(0,\alpha
\wedge2)$}{Proof of Theorem~1.1 for $r\in(0,\alpha
\wedge2)$}}\label{ss:0<r<alpha,2}

In this subsection, we prove Theorem~\ref{theorem:main1} for
$r\in(0,\alpha\wedge2)$. We will discuss the case for $\alpha\ne2$
and that
for $\alpha=2$ simultaneously, until we arrive at the point where we require
separate approaches.

First, we recall (\ref{eq:representation}) and split $\int_0^\infty$ into
$\int_0^U$ and $\int_U^\infty$ for a given $U>0$. Using
(\ref{eq:derivative-example}) for the former integral [as in
(\ref{eq:derivative-appl})] and (\ref{eq:laplace}) for
the latter, we obtain
%
\begin{eqnarray}\label{eq:small-r:dec1}
&&\sum_{t=0}^\infty m^t\sum_{x\in{\mathbb Z}^d}|x_1|^r\varphi
_t(x)\nonumber
\\
&&\qquad=\frac
{\hat\varphi(0)}
{K_r}\int_0^U\frac{\mathrm{d}u}{u^{1+r}} \frac{\bar\Delta_{\vec
u}\hat J(0)}{1-
\hat J(\vec u)}+\frac1{K_r}\int_0^U\frac{\mathrm{d}u}{u^{1+r}}
\frac{\bar\Delta_{\vec u}\hat I(0)}{1-\hat J(\vec u)}
\\
&&\qquad\quad{}+\frac1{K_r}\int_U^\infty\frac{\mathrm{d}u}{u^{1+r}}\sum
_{t=0}^\infty m^t
\sum_{x\in{\mathbb Z}^d}\bigl(1-\cos(ux_1)\bigr)\varphi_t(x).\nonumber
\end{eqnarray}
We note that, by (\ref{eq:Ihat}) and (\ref{eq:property3}),
$\bar\Delta_k\hat I(0)\equiv0$ for random walk and self-avoiding
walk and
$\bar\Delta_k\hat I(0)=O(\bar\Delta_k\hat D(0))$ uniformly in $m\le
m_{\rm c}$ for
oriented percolation. Since $\bar\Delta_k\hat J(0)$ is also
$O(\bar\Delta_k\hat D(0))$ uniformly in $m\le m_{\rm c}$ [see (\ref
{eq:Jhat-laplace})],
the integrals in the first two terms of (\ref{eq:small-r:dec1}) are of
the same
order and therefore the first term dominates the second term as
$m\uparrow m_{\rm c}$, due to the extra factor $\hat\varphi(0)$,
which exhibits
%
\begin{eqnarray}\label{eq:chi-asy}
\qquad \hat\varphi(0)&=&\frac{\hat I(0)}{\hat J_{\mathrm c}(0)-\hat
J(0)}=\frac
{\hat I_{\mathrm c}(0)+O
(1-{m}/{m_{\rm c}})}{m_{\rm c}\,\partial_m\hat J_{\mathrm c}(0)
(1-
{m}/{m_{\rm c}})+O((1-{m}/{m_{\rm c}})^{1+\epsilon})}\nonumber
\\[-8pt]\\[-8pt]
&=&\frac{C_{\mathrm I}}{1-{m}/{m_{\rm c}}}+O\biggl(\biggl(1-\frac
{m}{m_{\rm
c}}\biggr)^{-1+\epsilon}\biggr),\nonumber
\end{eqnarray}
where the first equality is due to (\ref{eq:lace-exp-solution}) and
(\ref{eq:critpt}), and the second equality is due to (\ref
{eq:property2}) and
(\ref{eq:Jhat-diff}). These estimates are valid independently of $r$
and thus
used in the later sections as well. By the fact that $0\le1-\cos
(ux_1)\le2$, the last term
in (\ref{eq:small-r:dec1}) obeys
\begin{eqnarray}\label{eq:small-r:dec4}
0&\le&\frac1{K_r}\int_U^\infty\frac{\mathrm{d}u}{u^{1+r}}\sum
_{t=0}^\infty m^t
\sum_{x\in{\mathbb Z}^d}\bigl(1-\cos(ux_1)\bigr)\varphi_t(x)\nonumber
\\[-8pt]\\[-8pt]
&\le&
\frac
{2\hat\varphi(0)}
{K_r}\int_U^\infty\frac{\mathrm{d}u}{u^{1+r}}
=\frac{2\hat\varphi(0)}{K_rr}U^{-r}.\nonumber
\end{eqnarray}
We will choose $U$ to be relatively small so as to make the first term in
(\ref{eq:small-r:dec1}) dominant.

Next, we investigate the integral part of the first term in
(\ref{eq:small-r:dec1}),
%
\begin{eqnarray}\label{eq:small-r:dec3}
\int_0^U\frac{\mathrm{d}u}{u^{1+r}} \frac{\bar\Delta_{\vec u}\hat
J(0)}{1-\hat J
(\vec u)}=\int_0^U\frac{\mathrm{d}u}{u^{1+r}} \frac{\bar\Delta
_{\vec
u}\hat J(0)}
{\hat J_{\mathrm c}(0)-\hat J(0)+\bar\Delta_{\vec u}\hat J(0)},
\end{eqnarray}
where we have used (\ref{eq:critpt}). By~(\ref{eq:D-asympt}) and
(\ref{eq:M-source}),
we have that, for small $u$,
\begin{eqnarray*}
\bar\Delta_{\vec u}\hat J(0)=\frac{\bar\Delta_{\vec u}\hat
J(0)}{\bar\Delta_{\vec u}
\hat D(0)} \bar\Delta_{\vec u}\hat D(0)=
\cases{
Mv_\alpha u^{\alpha\wedge2}+O(u^{\alpha\wedge2+\epsilon}),
&\quad$\alpha\ne2$,\cr
Mv_2u^2\log\dfrac1u+O(u^2),&\quad$\alpha=2$
}
\end{eqnarray*}
for some $\epsilon>0$, where the error terms are uniform in $m\le
m_{\rm c}$. Let
%
\begin{eqnarray}\label{eq:mu-def}
\mu=\frac{\hat J_{\mathrm c}(0)-\hat J(0)}{Mv_\alpha}.
\end{eqnarray}
Then,
\begin{eqnarray}\label{eq:small-r:dec5}
&&\frac{\bar\Delta_{\vec u}\hat J(0)}{\hat J_{\mathrm c}(0)-\hat
J(0)+\bar\Delta
_{\vec u}\hat
J(0)}\nonumber
\\[-8pt]\\[-8pt]
&&\qquad=
\cases{
\displaystyle\dfrac{u^{\alpha\wedge2}}{\mu+u^{\alpha\wedge
2}}+\frac
{O(u^{\alpha\wedge2
+\epsilon})}{\mu+u^{\alpha\wedge2}},&\quad$\alpha\ne2$,
\cr
\displaystyle\frac{u^2\log1/u}{\mu+u^2\log1/u}+\frac
{O(u^2)}{\mu+u^2\log1/
u},&\quad$\alpha=2$.
}\nonumber
\end{eqnarray}

We now investigate the integral (\ref{eq:small-r:dec3}) for $\alpha
\ne
2$ and
$\alpha=2$ separately, using (\ref{eq:small-r:dec5}) and the following
proposition.

\begin{proposition}\label{proposition:r<alpha&2}
Under the same conditions as in Theorem~\ref{theorem:main1},
%
\begin{eqnarray}\label{eq:M-asy}
M&=&M_{\mathrm c}+O\biggl(\biggl(1-\frac{m}{m_{\rm c}}\biggr)^\epsilon
\biggr),
\\\label{eq:mu-asy}
\mu&=&\frac{1-{m}/{m_{\rm c}}}{C_{\mathrm{II}}v_\alpha}+O
\biggl(\biggl(1-\frac
{m}{m_{\rm c}}\biggr)^{1+
\epsilon}\biggr)
\end{eqnarray}
for some $\epsilon>0$, where $M_{\mathrm c}=M(m_{\rm c})$.
\end{proposition}

The proof is deferred to Section~\ref{s:key-prop}. We note that these
estimates are trivial for random walk.

\subsubsection{Proof for $\alpha\ne2$}

\begin{figure}

\includegraphics{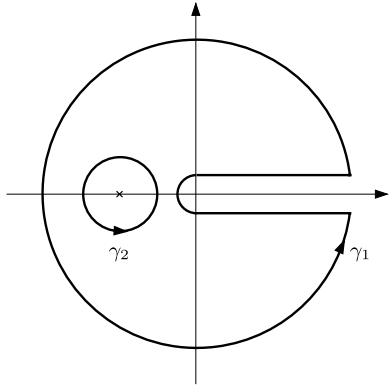}

\caption{The contours $\gamma_1$ and $\gamma_2$
in the complex plane.}\label{fig:contour1}
\end{figure}

We assume that $\epsilon<r$, without loss of generality. By~(\ref
{eq:small-r:dec3})
and (\ref{eq:small-r:dec5}) for $\alpha\ne2$, we have that, for
small $U$,
\begin{eqnarray*}
\int_0^U\frac{\mathrm{d}u}{u^{1+r}} \frac{\bar\Delta_{\vec u}\hat
J(0)}{1-\hat
J(\vec u)}&=&\int_0^U\frac{\mathrm{d}u}u\biggl(\frac{u^{\alpha\wedge
2-r}}{\mu
+u^{\alpha\wedge2}}+\frac{O(u^{\alpha\wedge2-r+\epsilon})}{\mu
+u^{\alpha
\wedge2}}\biggr)
\\
&=&\int_0^\infty\frac{\mathrm{d}u}u \frac{u^{\alpha\wedge2-r}}{\mu
+u^{\alpha
\wedge2}}-\int_U^\infty\frac{\mathrm{d}u}u \frac{u^{\alpha\wedge
2-r}}{\mu
+u^{\alpha\wedge2}}
\\
&&{}+\int_0^U\frac{\mathrm{d}u}u \frac{O(u^{\alpha\wedge
2+\epsilon-r})}
{\mu+u^{\alpha\wedge2}}\bigl({\mathbh1}_{\{\mu\ge u^{\alpha\wedge
2}\}}+{\mathbh1}_{\{\mu<u^{\alpha \wedge2}\}}\bigr)
\\
&=&\int_0^\infty\frac{\mathrm{d}u}u \frac{u^{\alpha\wedge2-r}}{\mu
+u^{\alpha
\wedge2}}+O(U^{-r})+O\bigl(\mu^{-{(r-\epsilon)/(\alpha\wedge2)}}\bigr).
\end{eqnarray*}
Let $U=\mu^{{(1-\epsilon/r)/(\alpha\wedge2)}}$, which is indeed
small as
$m\uparrow m_{\rm c}$, due to Proposition~\ref{proposition:r<alpha&2}.
By the change
of variables $u^{\alpha\wedge2}=\mu z$, we obtain
%
\begin{eqnarray}\label{eq:standard-cauchy}
&&\int_0^{\mu^{{(1-\epsilon/r)/(\alpha\wedge2)}}}\frac{\mathrm
{d}u}{u^{1+r}}
\frac{\bar\Delta_{\vec u}\hat J(0)}{1-\hat J(\vec u)}\nonumber
\\
&&\qquad =\int
_0^\infty
\frac{\mathrm{d}u}u \frac{u^{\alpha\wedge2-r}}{\mu+u^{\alpha
\wedge2}}
+O\bigl(\mu^{-{(r-\epsilon)/(\alpha\wedge2)}}\bigr)
\\
&&\qquad=\frac{\mu^{-{r/(\alpha\wedge2)}}}{\alpha\wedge2}\int
_0^\infty\frac{\mathrm
{d}z}z \frac{z^{1-{r/(\alpha\wedge2)}}}{1+z}+O\bigl(\mu^{-
{(r-\epsilon)/(\alpha\wedge2)}}\bigr).\nonumber
\end{eqnarray}
However, by the standard Cauchy integral formula, for $\beta\in(0,1)$,
%
\begin{eqnarray}\label{eq:standard-cauchy1}
\oint_{\gamma_1}\frac{\mathrm{d}z}z \frac{z^{1-\beta}}{1+z}=\oint
_{\gamma_2}
\frac{\mathrm{d}z}z \frac{z^{1-\beta}}{1+z}=2\pi i(-1)^{-\beta
}=2\pi i
e^{-\pi i\beta},
\end{eqnarray}
where, as depicted in Figure~\ref{fig:contour1}, the contour $\gamma_1$
consists of two line segments, an arc of the circle with smaller radius
$\delta\in(0,1)$ and an arc of the circle with larger radius $R\in
(1,\infty)$, and
the contour $\gamma_2$ is the circle centered at $-1$ with radius smaller
than~1. On the other hand, by taking $\delta\to0$ and $R\to\infty$,
we obtain
\begin{eqnarray*}
\mathop{\lim_{R\to\infty}}_{\delta\to0}\oint_{\gamma_1}\frac
{\mathrm{d}z}z
\frac{z^{1-\beta}}{1+z}
=(1-e^{-2\pi i\beta})\int_0^\infty\frac{\mathrm{d}z}z \frac
{z^{1-\beta}}{1+z}.
\end{eqnarray*}
Therefore,
\begin{eqnarray*}
\int_0^\infty\frac{\mathrm{d}z}z \frac{z^{1-\beta}}{1+z}=\frac
{2\pi i e^{-\pi
i\beta}}{1-e^{-2\pi i\beta}}=\frac\pi{\sin(\beta\pi)},
\end{eqnarray*}
which implies that
\begin{eqnarray*}
\int_0^{\mu^{{(1-\epsilon/r)/(\alpha\wedge2)}}}\frac{\mathrm
{d}u}{u^{1+r}}
\frac{\bar\Delta_{\vec u}\hat J(0)}{1-\hat J(\vec u)}&=&\frac\pi
{(\alpha\wedge2)
\sin(r\pi/(\alpha\wedge2))} \mu^{-{r/(\alpha\wedge 2)}}
\\
&&{}+O\bigl(\mu^{-
{(r-\epsilon)/(\alpha\wedge2)}}\bigr).
\end{eqnarray*}

Finally, by substituting (\ref{eq:standard-cauchy}) back into (\ref
{eq:small-r:dec1})
and using (\ref{eq:chi-asy}) and (\ref{eq:mu-asy}), we conclude that there
is an
$\epsilon'\in(0,1)$ such that
%
\begin{eqnarray}\label{eq:pr4r<2}
\qquad\ \ \sum_{t=0}^\infty m^t\sum_{x\in{\mathbb Z}^d}|x_1|^r\varphi
_t(x)&=&\frac{\pi
K_r^{-1}}
{(\alpha\wedge2)\sin(r\pi/(\alpha\wedge2))} \frac
{C_{\mathrm I}
(C_{\mathrm{II}}v_\alpha)
^{{r/(\alpha\wedge2)}}}{(1-{m/m_{\rm c}})^{1+
{r/(\alpha\wedge2)}}}\nonumber
\\[-8pt]\\[-8pt]
&&{}+O\biggl(\biggl(1-\frac{m}{m_{\rm c}}\biggr)^{-1-{r/(\alpha\wedge
2)}+\epsilon'}\biggr).\nonumber
\end{eqnarray}
However, since (see Appendix~\hyperref[appendix:K]{A.2})
%
\begin{eqnarray}\label{eq:Kr-evaluation}
\pi K_r^{-1}=2\Gamma(r+1)\sin\frac{r\pi}2,
\end{eqnarray}
this completes the proof of Theorem~\ref{theorem:main1} for
$0<r<\alpha\wedge2$ with $\alpha\ne2$.

\begin{Remark*}
Although the proportionality constant
$(2\sin\frac{r\pi}{\alpha\vee2})/((\alpha\wedge2)\sin\frac
{r\pi}\alpha)$ in
(\ref{eq:main1}) looks slightly different from the constant
$(2\sin\frac{r\pi}2)/((\alpha\wedge2)\sin\frac{r\pi}{\alpha
\wedge2})$ derived
from (\ref{eq:pr4r<2}) and (\ref{eq:Kr-evaluation}), they are equal when
$0<r<\alpha\wedge2$. The reason why we have adopted the former in the
main theorem is due to its applicability to larger values of $r<\alpha
$, which
the latter lacks (e.g., take $r=3<\alpha$).
\end{Remark*}

\subsubsection{Proof for $\alpha=2$}
The proof for $\alpha=2$ is slightly more involved than the above
proof for
$\alpha\ne2$, due to the log corrections in (\ref{eq:small-r:dec5}). By
(\ref{eq:small-r:dec3}) and (\ref{eq:small-r:dec5}) for $\alpha=2$, we
have that,
for small $U$,
\begin{eqnarray*}
\int_0^U\frac{\mathrm{d}u}{u^{1+r}} \frac{\bar\Delta_{\vec u}\hat
J(0)}{1-\hat J
(\vec u)}&=&\int_0^U\frac{\mathrm{d}u}u\biggl(\frac{u^{2-r}\log
1/u}{\mu+u^2
\log1/u}+\frac{O(u^{2-r})}{\mu+u^2\log1/u}\biggr)
\\
&=&\int_0^U\frac{\mathrm{d}u}u \frac{u^{2-r}\log1/u}{\mu
+u^2\log1/u}
+\frac{O(U^{2-r})}\mu,
\end{eqnarray*}
where we have obtained the error term by simply ignoring $u^2\log\frac1u>0$
in the denominator. Let $U=\sqrt\mu$, which is small as $m\uparrow
m_{\rm c}$, as
required, due to Proposition~\ref{proposition:r<alpha&2}. By the
change of
variables $u^2\log\frac1u=\mu z$, we obtain
\begin{eqnarray*}
\int_0^{\sqrt\mu}\frac{\mathrm{d}u}{u^{1+r}} \frac{\bar\Delta
_{\vec
u}\hat J(0)}
{1-\hat J(\vec u)}&=&\int_0^{\sqrt\mu}\frac{\mathrm{d}u}u \frac
{u^{2-r}\log
1/u}{\mu+u^2\log1/u}+O(\mu^{-r/2})
\\
&=&\frac{\mu^{-r/2}}2\int_0^{\log1/{\sqrt\mu}}\frac{\mathrm
{d}z}z \frac{z^{1
-{r}/2}(\log1/{u(z)})^{{r}/2}}{1+z}+O(\mu^{-r/2}).
\end{eqnarray*}
Note that, by taking the logarithm of $u^2\log1/u=\mu z$ and
using the
monotonicity of $(\log\log1/u)/\log1/u$ in $0<u<\sqrt\mu
\ll1$, we have
\begin{eqnarray*}
\log\frac1{u(z)}=\biggl(1+O\biggl(\frac{\log\log1/{\sqrt\mu
}}{\log1/
\mu}\biggr)\biggr)\log\frac1{\sqrt{\mu z}}.
\end{eqnarray*}
Therefore,
\begin{eqnarray*}
&&\int_0^{\sqrt\mu}\frac{\mathrm{d}u}{u^{1+r}} \frac{\bar\Delta
_{\vec u}\hat J(0)}
{1-\hat J(\vec u)}
\\
&&\qquad=\frac{\mu^{-r/2}}2\biggl(1+O\biggl(\frac{\log\log1/{\sqrt
\mu}}{\log1/
\mu}\biggr)\biggr)\int_0^{\log1/{\sqrt\mu}}\frac{\mathrm
{d}z}z \frac{z^{1-
{r}/2}(\log1/{\sqrt{\mu z}})^{{r}/2}}{1+z}
\\
&&{}\qquad\quad +O(\mu^{-r/2}).
\end{eqnarray*}
Suppose that $\log\frac1{\sqrt\mu}\gg1$. Then, by the Cauchy
integral formula (see
Figure~\ref{fig:contour1}),
\begin{eqnarray*}
\oint_{\gamma_1}\frac{\mathrm{d}z}z \frac{z^{1-{r}/2}(\log
1/{\sqrt{\mu
z}})^{{r}/2}}{1+z}&=&\oint_{\gamma_2}\frac{\mathrm{d}z}z \frac
{z^{1-{r}/
2}(\log1/{\sqrt{\mu z}})^{{r}/2}}{1+z}
\\
&=&2\pi i e^{-\pi ir/2}\biggl(\log\frac1{\sqrt\mu}-\frac{\pi
i}2\biggr)^{r/2}
\\
&=&2\pi i e^{-\pi ir/2}\biggl(\log\frac1{\sqrt\mu}\biggr)^{r/2}+O(1),
\end{eqnarray*}
where, as in (\ref{eq:standard-cauchy1}), the contour $\gamma_2$ is
the circle
at $-1$ with radius smaller than 1, while the contour $\gamma_1$
contains an
arc of the circle with radius $\delta\in(0,1)$ and an arc of the
circle with
radius $R\equiv\log\frac1{\sqrt\mu}$. On the other hand, by taking
$\delta\to0$, we obtain
\begin{eqnarray*}
&&\lim_{\delta\to0}\oint_{\gamma_1}\frac{\mathrm{d}z}z \frac
{z^{1-{r}/2}(\log
1/{\sqrt{\mu z}})^{{r}/2}}{1+z}
\\
&&\qquad =(1-e^{-\pi ir})\int_0^{\log1/
{\sqrt\mu}}\frac{\mathrm{d}z}z \frac{z^{1-{r}/2}(\log
1/{\sqrt{\mu
z}})^{{r}/2}}{1+z}+O(1),
\end{eqnarray*}
where the error term is independent of $\mu$. Therefore,
\begin{eqnarray*}
\int_0^{\log1/{\sqrt\mu}}\frac{\mathrm{d}z}z \frac
{z^{1-{r}/2}(\log
1/{\sqrt{\mu z}})^{{r}/2}}{1+z}&=&\frac{2\pi i e^{-\pi ir/2}}{1-
e^{-\pi ir}}\biggl(\log\frac{1}{\sqrt\mu}\biggr)^{r/2}+O(1)
\\
&=&\frac\pi{\sin{(r\pi/2)}}\biggl(\log\frac1{\sqrt\mu}
\biggr)^{r/2}+O(1),
\end{eqnarray*}
which implies
that
%
\begin{eqnarray}\label{eq:generalized-cauchy}
\qquad \int_0^{\sqrt\mu}\frac{\mathrm{d}u}{u^{1+r}} \frac{\bar\Delta
_{\vec
u}\hat
J(0)}{1-\hat J(\vec u)}=\frac\pi{2\sin({r\pi}/2)} \mu
^{-r/2}\biggl(\log\frac1
{\sqrt\mu}\biggr)^{r/2}+O(\mu^{-r/2}),
\end{eqnarray}
where we have used
\begin{eqnarray*}
O\biggl(\frac{\log\log1/{\sqrt\mu}}{\log1/\mu}
\biggr)\biggl(\log\frac1
{\sqrt\mu}\biggr)^{r/2}=o(1)\qquad[\because r<2].
\end{eqnarray*}

Finally, by substituting (\ref{eq:generalized-cauchy}) back into
(\ref{eq:small-r:dec1}) and using (\ref{eq:small-r:dec4}) with
$U=\sqrt\mu$,
we obtain
\begin{eqnarray*}
\sum_{t=0}^\infty m^t\sum_{x\in{\mathbb Z}^d}|x_1|^r\varphi
_t(x)&=&\bigl(\hat
\varphi(0)
+O(1)\bigr)\frac{\pi K_r^{-1}}{2\sin({r\pi}/2)} \mu^{-r/2}
\biggl(\log\frac1
{\sqrt\mu}\biggr)^{r/2}
\\
&&{}+\hat\varphi(0) O(\mu^{-r/2}).
\end{eqnarray*}
Combining this with (\ref{eq:chi-asy}), (\ref{eq:mu-asy}) and (\ref
{eq:Kr-evaluation})
yields (\ref{eq:main1}) for $\alpha=2$. This completes the proof of
Theorem~\ref{theorem:main1} for $0<r<\alpha=2$.

\subsection{\texorpdfstring{Proof of Theorem~\protect\ref{theorem:main1} for $r=2j<\alpha$
[$j\in\mathbb{N}$]}{Proof of Theorem~1.1 for $r=2j<\alpha$
[$j\in\mathbb{N}$]}}\label{ss:r=2j}

In this subsection, we prove Theorem~\ref{theorem:main1} for positive even
integers $r=2j<\alpha$. First, we recall (\ref{eq:representation})
for $r=2j$:
\begin{eqnarray*}
\sum_{t=0}^\infty m^t\sum_{x\in{\mathbb Z}^d}x_1^{2j}\varphi_t(x)
&=(-1)^j\nabla_{1}^{2j}\hat\varphi(0).
\end{eqnarray*}
Differentiating (\ref{eq:lace-exp-fourier}) and using the ${\mathbb
Z}^d$-symmetry
of the
models [so that
$\nabla_{1}^n\hat J(0)$ and $\nabla_{1}^n\hat\varphi(0)$ are both zero
when $n$ is
odd], we have
\begin{eqnarray*}
\nabla_1^{2j}\hat\varphi(0)=\nabla_{1}^{2j}\hat I(0)+\hat J(0)
\nabla_{1}^{2j}\hat
\varphi(0)+\sum_{l=1}^j\pmatrix{2j\cr 2l}\nabla_{1}^{2l}\hat J(0)
\nabla_{1}^{2(j-l)}
\hat\varphi(0).
\end{eqnarray*}
Solving this equation for $\nabla_1^{2j}\hat\varphi(0)$ and using
(\ref{eq:lace-exp-solution}) for $k=0$, we obtain
%
\begin{eqnarray}\label{eq:nablone2j-lace}
\nabla_1^{2j}\hat\varphi(0)=\frac{\hat\varphi(0)}{\hat I(0)}
\Biggl(\nabla_{1}^{2j}
\hat I(0)+\sum_{l=1}^j\pmatrix{2j\cr 2l}\nabla_{1}^{2l}\hat J(0) \nabla_{1}
^{2(j-l)}
\hat\varphi(0)\Biggr).
\end{eqnarray}

To identify the dominant term of the right-hand side, we use the following
proposition.

\begin{proposition}\label{proposition:2<r<alpha}
Let $\alpha>2$ and $\langle\alpha\rangle=\max\{j\in\mathbb
{N}\dvtx j<\alpha\}
$ (note that
$\langle\alpha\rangle=\alpha-1$ if $\alpha\ge3$ is an integer).
Under the
same conditions as in Theorem~\ref{theorem:main1},\vspace*{-2pt}
%
\begin{eqnarray}\label{eq:2<r<alpha-IJbds}
\qquad \quad \left.\begin{array}{r}
\displaystyle\sum_{t=0}^\infty m^t\sum_{x\in{\mathbb Z}^d}|x_1|^\nu
|I_t(x)|\\
\displaystyle\sum_{t=0}^\infty m^t\sum_{x\in{\mathbb Z}^d}|x_1|^\nu|J_t(x)|
\end{array}
\right\}\le
\cases{
O(1),&\quad$0\le\nu\le2$,
\cr
O\biggl(\biggl(1-\dfrac{m}{m_{\rm c}}\biggr)^{1-\nu/2+\epsilon}\biggr),
&\quad$2<\nu<\alpha$
}\hspace*{-15pt}\vspace*{-2pt}
\end{eqnarray}
for some $\epsilon>0$. Moreover,\vspace*{-2pt}
\begin{eqnarray}\label{eq:2<r<alpha-varphibd}
\quad  &&|\nabla_{1}^n\hat\varphi(k,me^{i\theta})|\nonumber
\\[-9pt]\\[-9pt]
 &&\qquad \le
\cases{
O\biggl(\biggl(1-\dfrac{m}{m_{\rm c}}+|\theta|+|k|^2\biggr)^{-1-{n}/2}
\biggr),
&\quad$n=0,1,2$,
\cr
\dfrac{O((1-{m}/{m_{\rm c}})^{1-{n}/2}
)}{(1-{m}/{m_{\rm c}}
+|\theta|+|k|^2)^2},&\quad$n=3,\dots,\langle\alpha\rangle$,
}\nonumber\vspace*{-2pt}
\end{eqnarray}
where the $O((1-\frac{m}{m_{\rm c}})^{1-{n}/2})$ term is uniform in
$(k,\theta)\in[-\pi,\pi]^{d+1}$.
\end{proposition}

We will use this proposition again in the next subsection to prove
Theorem~\ref{theorem:main1} for the remaining case: $r=2j+q$, where
$j\in\mathbb{N}$
and $q\in(0,2)$. The proof of Proposition~\ref{proposition:2<r<alpha} is
deferred to Section~\ref{s:key-prop}. Note that (\ref{eq:2<r<alpha-IJbds})
is trivial for random walk.

Now we resume the proof of Theorem~\ref{theorem:main1} for $r=2j$.
Notice that\vspace*{-2pt}
%
\begin{eqnarray}\label{eq:nablone2lJ}
|\nabla_{1}^{2l}\hat J(0)|\le\sum_{t=0}^\infty m^t\sum_{x\in
{\mathbb Z}^d}|x_1|^{2l}
|J_t(x)|\vspace*{-2pt}
\end{eqnarray}
and that a similar bound holds for $I$. By
(\ref{eq:nablone2j-lace})--(\ref{eq:2<r<alpha-varphibd}), we have the
recursion\vspace*{-2pt}
\begin{eqnarray*}
\nabla_{1}^{2j}\hat\varphi(0)&=&\frac{\hat\varphi(0)}{\hat
I(0)}\Biggl(\nabla_{1}^{2
j}\hat I(0)+\pmatrix{2j\cr 2}\nabla_{1}^2\hat J(0) \nabla
_{1}^{2(j-1)}\hat
\varphi(0)
\\[-2pt]
&&{}\hspace*{50pt}+\sum_{l=2}^j\pmatrix{2j\cr 2l}\nabla_{1}^{2l}\hat J(0)
\nabla_{1}^{2(j-l)}
\hat\varphi(0)\Biggr)
\\[-2pt]
&=&\pmatrix{2j\cr 2}\frac{\nabla_{1}^2\hat J(0)}{\hat I(0)} \hat\varphi
(0) \nabla_{1}^{2
(j-1)}\hat\varphi(0)+O\biggl(\biggl(1-\frac{m}{m_{\rm c}}\biggr)^{-1-j+\epsilon
}\biggr),\vspace*{-2pt}
\end{eqnarray*}
where the first term is $O((1-\frac{m}{m_{\rm c}})^{-1-j})$, which is
dominant as
$m\uparrow m_{\rm c}$. Repeated use of this recursion then yields
\begin{eqnarray*}
\nabla_{1}^{2j}\hat\varphi(0)&=&\pmatrix{2j\cr2}\pmatrix{2(j-1)\cr 2}
\biggl(\frac{\nabla_{1}^2
\hat J(0)}{\hat I(0)} \hat\varphi(0)\biggr)^2\nabla
_{1}^{2(j-2)}\hat
\varphi
(0)
\\[-2pt]
&&{}+O\biggl(\biggl(1-\frac{m}{m_{\rm c}}\biggr)^{-1-j+\epsilon}\biggr)
\\
&\vdots&
\\[-2pt]
&=&\prod_{l=2}^j\pmatrix{2l\cr 2}\biggl(\frac{\nabla_{1}^2\hat J(0)}{\hat
I(0)} \hat
\varphi(0)\biggr)^{j-1}\nabla_{1}^2\hat\varphi(0)+O\biggl(\biggl(1-\frac
{m}{m_{\rm c}}\biggr)^{-1-j
+\epsilon}\biggr)
\\[-2pt]
&=&\frac{(2j)!}{2^j}\biggl(\frac{\nabla_{1}^2\hat J(0)}{\hat
I(0)}
\biggr)^j\hat
\varphi(0)^{j+1}+O\biggl(\biggl(1-\frac{m}{m_{\rm c}}\biggr)^{-1-j+\epsilon}\biggr).
\end{eqnarray*}
However, by comparing (\ref{eq:nablone}) and (\ref{eq:laplace}), and using
(\ref{eq:M-asy}),
we have
\begin{eqnarray*}
\nabla_{1}^2\hat J(0)&=&-2v_\alpha\lim_{k\to0}\frac{\bar\Delta
_k\hat
J(0)}{\bar\Delta_k
\hat D(0)}=-2v_\alpha M
\\[-2pt]
&=&-2v_\alpha M_{\mathrm c}+O\biggl(\biggl(1-\frac{m}{m_{\rm c}}\biggr)^\epsilon
\biggr).
\end{eqnarray*}
Recall that $\hat I(0)=\hat I_{\mathrm c}(0)+O(1-\frac{m}{m_{\rm c}})$
[cf.~the numerator
in (\ref{eq:chi-asy})]. Therefore,
%
\begin{eqnarray}\label{eq:ratio}
\quad \frac{\nabla_{1}^2\hat J(0)}{\hat I(0)}&=&-2v_\alpha\frac{M_{\mathrm
c}}{\hat
I_{\mathrm c}(0)}
+O\biggl(\biggl(1-\frac{m}{m_{\rm c}}\biggr)^\epsilon\biggr)\nonumber
\\[-9pt]\\[-9pt]
&=&-2v_\alpha\frac{C_{\mathrm{II}}}{C_{\mathrm I}}+O\biggl(\biggl(1-\frac
{m}{m_{\rm
c}}\biggr)^\epsilon\biggr)\qquad
[\because(\ref{eq:CICII})\mbox{ and }(\ref{eq:M-source})],\nonumber
\end{eqnarray}
hence
\begin{eqnarray*}
\nabla_{1}^{2j}\hat\varphi(0)&=&\frac{(2j)!}{2^j}\biggl(-2v_\alpha
\frac{C_{\mathrm{II}}}{C_{\mathrm I}}
\biggr)^j\biggl(\frac{C_{\mathrm I}}{1-{m/m_{\rm c}}}
\biggr)^{j+1}+O\biggl(\biggl(1-\frac{m}
{m_{\rm c}}\biggr)^{-1-j+\epsilon}\biggr)
\\[-2pt]
&=&\Gamma(2j+1) \frac{C_{\mathrm I}(-C_{\mathrm{II}}v_\alpha
)^j}{(1-{m/m_{\rm
c}})^{j+1}}+O\biggl(
\biggl(1-\frac{m}{m_{\rm c}}\biggr)^{-1-j+\epsilon}\biggr).
\end{eqnarray*}

This completes the proof of Theorem~\ref{theorem:main1} for positive even
integers $r=2j<\alpha$.

\subsection{\texorpdfstring{Proof of Theorem~\protect\ref{theorem:main1} for $r=2j+q<\alpha$
[$j\in\mathbb{N}$, $q\in(0,2)$]}{Proof of Theorem~1.1 for $r=2j+q<\alpha$
[$j\in\mathbb{N}$, $q\in(0,2)$]}}\label{ss:r=2j+q}

In this subsection, we prove Theorem~\ref{theorem:main1} for the other values
of $r<\alpha$: $r=2j+q$ with $j\in\mathbb{N}$ and $q\in(0,2)$.
First, we recall
(\ref{eq:representation}):
%
\begin{eqnarray}\label{eq:pr42j+q:1}
\sum_{t=0}^\infty m^t\sum_{x\in{\mathbb Z}^d}|x_1|^r\varphi
_t(x)=\frac{(-1)^j}{K_q}
\int_0^\infty\frac{\mathrm{d}u}{u^{1+q}} \bar\Delta_{\vec u}\nabla_{1}
^{2j}\hat
\varphi(0),
\end{eqnarray}
where, by (\ref{eq:lace-exp-fourier}),
\begin{eqnarray*}
\bar\Delta_{\vec u}\nabla_{1}^{2j}\hat\varphi(0)&=&\nabla
_{1}^{2j}\hat
\varphi(0)
-\nabla_{1}^{2j}\hat\varphi(\vec u)
\\[-2pt]
&=&\nabla_{1}^{2j}\hat I(0)-\nabla_{1}^{2j}\hat I(\vec u)+\sum
_{n=0}^{2j}\pmatrix{2j\cr
n}\bigl(\nabla_{1}^n\hat J(0) \nabla_{1}^{2j-n}\hat\varphi(0)
\\[-2pt]
&&{}\hspace*{150pt}-\nabla_{1}^n\hat J(\vec u) \nabla_{1}^{2j-n}\hat
\varphi
(\vec u)\bigr)
\\
&=&\bar\Delta_{\vec u}\nabla_{1}^{2j}\hat I(0)+\sum_{n=0}^{2j}
\pmatrix{2j\cr n }\bigl(
\nabla_{1}^n\hat J(0)~\bar\Delta_{\vec u}\nabla_{1}^{2j-n}\hat
\varphi (0)
\\
&&{}\hspace*{112pt}+\nabla_{1}^{2j-n}\hat\varphi(\vec u) \bar\Delta_{\vec
u}\nabla_{1}^n
\hat J(0)\bigr).
\end{eqnarray*}
Solving this equation for $\bar\Delta_{\vec u}\nabla_{1}^{2j}\hat
\varphi
(0)$ and
using (\ref{eq:lace-exp-solution}) for $k=0$ and $\nabla_{1}^n\hat J(0)=0$
for odd
$n$, we obtain
\begin{eqnarray*}
\bar\Delta_{\vec u}\nabla_{1}^{2j}\hat\varphi(0)&=&\frac{\hat
\varphi
(0)}{\hat I(0)}
\Biggl(\bar\Delta_{\vec u}\nabla_{1}^{2j}\hat I(0)+\sum
_{l=1}^j\pmatrix{2j\cr 2l}
\nabla_{1}^{2l}\hat J(0) \bar\Delta_{\vec u}\nabla_{1}^{2(j-l)}\hat
\varphi(0)
\\
&&{}\hspace*{86pt}+\sum_{n=0}^{2j}\pmatrix{2j\cr n}\nabla_{1}^{2j-n}\hat\varphi(\vec
u) \bar\Delta_{\vec u}
\nabla_{1}^n\hat J(0)\Biggr).
\end{eqnarray*}
Substituting this back into (\ref{eq:pr42j+q:1}) yields
\begin{eqnarray}\label{eq:pr42j+q:2}
&&\sum_{t=0}^\infty m^t\sum_{x\in{\mathbb Z}^d}|x_1|^r\varphi\nonumber
_t(x)
\\[-8pt]\\[-8pt]
&&\qquad=\frac{\hat
\varphi(0)}
{\hat I(0)}\Biggl(H^{(1)}+\sum_{l=1}^j\pmatrix{2j\cr 2l}H_{2l}^{(2)}+
\sum_{n=0}^{2j}\pmatrix{2j\cr n}H_n^{(3)}\Biggr),\nonumber
\end{eqnarray}
where
%
\begin{eqnarray}\label{eq:H1-def}
H^{(1)}&=&\frac{(-1)^j}{K_q}\int_0^\infty\frac{\mathrm
{d}u}{u^{1+q}}
\bar\Delta_{\vec u}\nabla_{1}^{2j}\hat I(0)\equiv\sum_{t=0}^\infty m^t
\sum_{x\in{\mathbb Z}^d}|x_1|^rI_t(x),
\\\label{eq:H2-def}
H_{2l}^{(2)}&=&\frac{(-1)^j}{K_q}\int_0^\infty\frac{\mathrm
{d}u}{u^{1+q}}
\nabla_{1}^{2l}\hat J(0) \bar\Delta_{\vec u}\nabla_{1}^{2(j-l)}\hat
\varphi(0)\nonumber
\\[-8pt]\\[-8pt]
&\equiv&(-1)^l\nabla_{1}^{2l}\hat J(0)\sum_{t=0}^\infty m^t\sum
_{x\in
{\mathbb Z}^d}|x_1|^{r
-2l}\varphi_t(x)\nonumber
\end{eqnarray}
and
%
\begin{eqnarray}\label{eq:H3-def}
H_n^{(3)}&=&\frac{(-1)^j}{K_q}\int_0^\infty\frac{\mathrm
{d}u}{u^{1+q}}
\nabla_{1}^{2j-n}\hat\varphi(\vec u)~\bar\Delta_{\vec u}\nabla
_{1}^n\hat
J(0)\nonumber
\\
&\equiv&\sum_{s,t=0}^\infty m^{t+s}\sum_{x,y\in{\mathbb Z}^d
}x_1^{2j-n}\varphi_t(x)
y_1^nJ_s(y)
\\
&&{}\times\frac1
{K_q}\int_0^\infty\frac{\mathrm{d}u}{u^{1+q}}\times
\cases{
\sin(ux_1) \sin(uy_1),&\quad\mbox{odd }$n$,
\cr
\cos(ux_1) \bigl(1-\cos(uy_1)\bigr),&\quad\mbox{even }$n$.
}
\nonumber
\end{eqnarray}

Next, we isolate error terms from (\ref{eq:pr42j+q:2}) using
Proposition~\ref{proposition:2<r<alpha}. First, by (\ref{eq:2<r<alpha-IJbds}),
we have
%
\begin{eqnarray}\label{eq:H1-bd}
\bigl|H^{(1)}\bigr|\le\sum_{t=0}^\infty m^t\sum_{x\in{\mathbb Z}^d
}|x_1|^r|I_t(x)|\le O\biggl(\biggl(1
-\frac{m}{m_{\rm c}}\biggr)^{-{r}/2+1+\epsilon}\biggr),
\end{eqnarray}
which gives rise to an error term.

Next, for $H_{2l}^{(2)}$, where $r-2l=2j+q-2l<2j+2-2l<\alpha$, we first
apply Jensen's inequality and then (\ref{eq:2<r<alpha-varphibd}) to obtain
%
\begin{eqnarray}\label{eq:Jensen}
\quad&&\sum_{t=0}^\infty m^t\sum_{x\in{\mathbb Z}^d}|x_1|^{r-2l}\varphi
_t(x)\nonumber
\\
&&\qquad \le
\Biggl(\frac1
{\hat\varphi(0)}\sum_{t=0}^\infty m^t\sum_{x\in{\mathbb Z}^d
}|x_1|^{2j+2-2l}\varphi_t
(x)\Biggr)^{{(r-2l)/(2j+2-2l)}}\hat\varphi(0)\nonumber
\\
&&\qquad=\biggl(\frac{|\nabla_{1}^{2j+2-2l}\hat\varphi(0)|}{\hat\varphi(0)}
\biggr)^{{(r-2l)/(2j+2-2l)}}\hat\varphi(0)
\\
&&\qquad\le O\biggl(\biggl(1-\frac{m}{m_{\rm c}}\biggr)^{-{(2j+2-2l)}/2}\biggr)^{
{(r-2l)/(2j+2-2l)}}
 O\biggl(\biggl(1-\frac{m}{m_{\rm c}}\biggr)^{-1}\biggr)\nonumber
 \\
&&\qquad =O\biggl(\biggl(1-\frac{m}{m_{\rm c}}\biggr)^{-1-{(r-2l)}/2}\biggr).\nonumber
\end{eqnarray}
Combining this with (\ref{eq:2<r<alpha-IJbds}) and (\ref
{eq:nablone2lJ}) yields
%
\begin{eqnarray}\label{eq:H2-bd}
\bigl|H_{2l}^{(2)}\bigr|\le
\cases{
O\biggl(\biggl(1-\dfrac{m}{m_{\rm c}}\biggr)^{-r/2}\biggr),&\quad$l=1$,
\cr
O\biggl(\biggl(1-\dfrac{m}{m_{\rm c}}\biggr)^{-{r}/2+\epsilon}\biggr),
&\quad$l=2,3,\dots,j$.
}
\end{eqnarray}

Finally, for $H_n^{(3)}$ with $n\ge2$ ($H_0^{(3)}$ and
$H_1^{(3)}$
will be investigated in detail later), we use
\begin{eqnarray}\label{eq:crudebd-odd}
\int_0^\infty\frac{\mathrm{d}u}{u^{1+q}} |\sin(ux_1) \sin(uy_1)|
&\le&\int_0^\infty\frac{\mathrm{d}u}{u^{1+q}} (|u^2x_1y_1|\wedge1)
\nonumber
\\[-8pt]\\[-8pt]
&=&O(|x_1y_1|^{q/2}),\nonumber
\\\label{eq:crudebd-even}
\int_0^\infty\frac{\mathrm{d}u}{u^{1+q}} \bigl|\cos(ux_1)
\bigl(1-\cos(uy_1)
\bigr)\bigr|
&\le&\int_0^\infty\frac{\mathrm{d}u}{u^{1+q}} \biggl(\frac{u^2y_1^2}2
\wedge2\biggr)\nonumber
\\[-8pt]\\[-8pt]
&=&O(|y_1|^q),\nonumber
\end{eqnarray}
which are due to the naive bounds $|\sin w|\le|w|\wedge1$, $|\cos
w|\le1$
and $|1-\cos w|\le\frac{w^2}2\wedge2$. By
(\ref{eq:crudebd-odd}) and (\ref{eq:crudebd-even}) and using Jensen's
inequality for
odd $n,$ as in (\ref{eq:Jensen}), we obtain
\begin{eqnarray*}
\bigl|H_n^{(3)}\bigr|
&\le&
\cases{
\biggl(\dfrac{|\nabla_{1}^{2j-n+1}\hat\varphi
(0)|}{\hat\varphi(0)}
\biggr)^{{(2j-n+{q}/2)/(2j-n+1)}}\hat\varphi(0)
\cr
\qquad {}\times \displaystyle\sum
_{s=0}^\infty
m^s\sum_{y\in{\mathbb Z}^d}|y_1|^{n+{q}/2}|J_s(y)|,&\quad\mbox{odd
}$n$,
\cr
\displaystyle|\nabla_{1}^{2j-n}\hat\varphi(0)|\sum_{s=0}^\infty
m^s\sum
_{y\in{\mathbb Z}^d}
|y_1|^{n+q}|J_s(y)|,&\quad\mbox{even }$n$.
}
\end{eqnarray*}
Then, by Proposition~\ref{proposition:2<r<alpha} and using $2j+q=r$,
we obtain
%
\begin{eqnarray}\label{eq:H3-bd}
\bigl|H_n^{(3)}\bigr|\le O\biggl(\biggl(1-\frac{m}{m_{\rm c}}\biggr)^{-
{r}/2+\epsilon}\biggr)\qquad
[n=2,3,\dots,2j].
\end{eqnarray}

Now, by (\ref{eq:pr42j+q:2}), (\ref{eq:H1-bd}), (\ref{eq:H2-bd}) and
(\ref{eq:H3-bd}), we arrive at
%
\begin{eqnarray}\label{eq:pr42j+q:3}
\quad \sum_{t=0}^\infty m^t\sum_{x\in{\mathbb Z}^d}|x_1|^r\varphi
_t(x)&=&\frac
{\hat\varphi(0)}
{\hat I(0)}\left(\pmatrix{2j\cr 2}H_2^{(2)}+H_0^{(3)}+\pmatrix
{2j\cr 1}H_1^{(3)}\right)\nonumber
\\[-8pt]\\[-8pt]
&&{}+O\biggl(\biggl(1-\frac{m}{m_{\rm c}}\biggr)^{-1-{r}/2+\epsilon}
\biggr).\nonumber
\end{eqnarray}

Finally, we reorganize the main term of (\ref{eq:pr42j+q:3}) and
complete the
proof of Theorem~\ref{theorem:main1}. First, we note that
\begin{eqnarray*}
\sin(ux_1) \sin(uy_1)&=&\frac{\cos(u(x_1-y_1))-\cos(u(x_1+y_1))}2
\\
&=&\frac{1-\cos(u(x_1+y_1))-(1-\cos(u(x_1-y_1))
)}2,\nonumber
\\
\cos(ux_1) \bigl(1-\cos(uy_1)\bigr)&=&\cos(ux_1)-\frac{\cos(u(x_1+y_1))
+\cos(u(x_1-y_1))}2
\\
&=&\frac{(1-\cos(u(x_1+y_1)))+(1-\cos(u(x_1-y_1))
)}2
\\
&&{}-\bigl(1-\cos(ux_1)\bigr).
\end{eqnarray*}
Then, by (\ref{eq:Kr}), we have the identities
\begin{eqnarray*}
\frac1{K_q}\int_0^\infty\frac{\mathrm{d}u}{u^{1+q}} \sin(ux_1)
\sin(uy_1)&=&
\frac{|x_1+y_1|^q-|x_1-y_1|^q}2,
\\
\frac1{K_q}\int_0^\infty\frac{\mathrm{d}u}{u^{1+q}} \cos(ux_1)
\bigl(1-\cos(u
y_1)\bigr)&=&\frac{|x_1+y_1|^q+|x_1-y_1|^q-2|x_1|^q}2.
\end{eqnarray*}
By these identities and the fact that $r=2j+q$, we obtain
\begin{eqnarray*}
\pmatrix{2j\cr 2}H_2^{(2)}+H_0^{(3)}+\pmatrix{2j\cr 1}H_1^{
(3)}=\sum_{s,t
=0}^\infty m^{t+s}\sum_{x,y\in{\mathbb Z}^d}x_1^{2j-2}\varphi_t(x)
J_s(y)
{\mathcal H}
(x_1,y_1),
\end{eqnarray*}
where
\begin{eqnarray*}
{\mathcal H}(x_1,y_1)&=&\pmatrix{2j\cr 2}|x_1|^qy_1^2+x_1^2\frac
{|x_1+y_1|^q+|x_1-y_1|^q
-2|x_1|^q}2
\\
&&{}+\pmatrix{2j\cr 1}x_1y_1\frac{|x_1+y_1|^q-|x_1-y_1|^q}2.
\end{eqnarray*}
In fact, due to the symmetry
${\mathcal H}(x_1,y_1)={\mathcal H}(x_1,-y_1)={\mathcal
H}(-x_1,y_1)={\mathcal H}(-x_1,\break -y_1)$ for any
$x_1,y_1\in{\mathbb Z}$, the above identity is equivalent to
\begin{eqnarray*}
&&\pmatrix{2j\cr 2}H_2^{(2)}+H_0^{(3)}+\pmatrix{2j\cr 1}H_1^{(3)}
\\
&&\qquad =4\sum_{s,t=0}^\infty m^{t+s}\mathop{\sum_{x,y\in{\mathbb
Z}^d}}_{
(x_1,y_1>0)}
x_1^{2j-2}\varphi_t(x) J_s(y) {\mathcal H}(x_1,y_1).
\end{eqnarray*}
Using the Taylor expansion of
$|x_1\pm y_1|^q\equiv x_1^q(1\pm\frac{y_1}{x_1})^q$ if $x_1>y_1>0$ and
the expansion of $|x_1\pm y_1|^q\equiv y_1^q(1\pm\frac{x_1}{y_1})^q$ if
$y_1>x_1>0$, we have
%
\begin{eqnarray}\label{eq:cH-identity}
&&{\mathcal H}(x_1,y_1)\nonumber
\\[-9pt]\\[-9pt]
&&\qquad=
\cases{
\displaystyle\left(\pmatrix{2j\cr 2}+\pmatrix{q\cr 2}+\pmatrix{2j\cr 1}q\right)x_1^qy_1^2
+O(x_1^{q-1}y_1^3),\cr \qquad x_1>y_1>0,
\cr
O(y_1^{2+q}),\qquad y_1\ge x_1>0.
}
\nonumber
\end{eqnarray}
Notice that
\begin{eqnarray*}
\pmatrix{2j\cr 2}+\pmatrix{q\cr 2}+\pmatrix{2j\cr 1}q&=&j(2j-1)+\frac{q}2(q-1)+2jq
\\
&=&\biggl(j+\frac{q}2\biggr)(2j+q)-j-\frac{q}2=\frac{r}2r-\frac{r}2
=\pmatrix{r\cr2}.
\end{eqnarray*}
We also notice that, as long as $q\in(0,1]$,
we have
\begin{eqnarray*}
x_1^{q-1}y_1^3=\biggl(\frac{y_1}{x_1}\biggr)^{1-q}y_1^{2+q}\le y_1^{2+q}
\qquad[x_1>y_1>0].
\end{eqnarray*}
Therefore, by Proposition~\ref{proposition:2<r<alpha}, we obtain that,
for $q\equiv r-2j\in(0,1]$,
%
\begin{eqnarray}\label{eq:pr42j+q:4.5}
&&\pmatrix{2j\cr 2}H_2^{(2)}+H_0^{(3)}+\pmatrix{2j\cr 1}H_1^{(3)}\nonumber
\\[-2pt]
&&\qquad=4\pmatrix{r\cr 2}\sum_{s,t=0}^\infty m^{t+s}\mathop{\sum_{x,y\in
{\mathbb Z}^d
}}_{ (x_1,y_1>
0)}x_1^{2j+q-2}\varphi_t(x) y_1^2J_s(y)\nonumber
\\[-2pt]
&&{}\qquad\quad+\sum_{s,t=0}^\infty m^{t+s}\sum_{x,y\in{\mathbb Z}^d
}x_1^{2j-2}\varphi_t(x)
O(|y_1|^{2+q})J_s(y)
\\[-2pt]
&&\qquad=\pmatrix{r\cr 2}(-\nabla_{1}^2\hat J(0))\sum_{t=0}^\infty
m^t\sum_{x\in{\mathbb Z}^d}
|x_1|^{r-2}\varphi_t(x)\nonumber
\\[-2pt]
&&\qquad{}\quad+\underbrace{|\nabla_{1}^{2j-2}\hat\varphi(0)|\sum
_{s=0}^\infty m^s\sum_{y
\in{\mathbb Z}^d}O(|y_1|^{2+q})J_s(y)}_{O((1-{m/m_{\rm
c}})^{-
{r}/2+\epsilon})}.\nonumber
\end{eqnarray}
For $q\in(1,2)$, we have to deal with the contribution from
$O(x_1^{q-1}y_1^3)$ in (\ref{eq:cH-identity}). However, by Jensen's inequality
and Proposition~\ref{proposition:2<r<alpha},
we have
\begin{eqnarray*}
&&\sum_{s,t=0}^\infty m^{t+s}\mathop{\sum_{x,y\in{\mathbb Z}^d
}}_{(x_1>y_1>0)} x_1^{2j
+q-3}\varphi_t(x) y_1^3|J_s(y)|
\\[-2pt]
&&\qquad \le\biggl(\frac{|\nabla_{1}^{2j-1}\hat\varphi(0)|}{\hat\varphi
(0)}\biggr)^{{(2
j+q-3)/(2j-1)}}\hat\varphi(0)\sum_{s=0}^\infty m^s\sum_{y\in{\mathbb Z}^d
}|y_1|^3|J_s(y)|
\\[-2pt]
&&\qquad \le O\biggl(\biggl(1-\frac{m}{m_{\rm c}}\biggr)^{-{r}/2+\epsilon}\biggr)
\end{eqnarray*}
and thus (\ref{eq:pr42j+q:4.5}) is valid for any $q\in(0,2)$.

Now, by substituting (\ref{eq:pr42j+q:4.5}) back into (\ref
{eq:pr42j+q:3}), we
obtain the recursion
\begin{eqnarray*}
\sum_{t=0}^\infty m^t\sum_{x\in{\mathbb Z}^d}|x_1|^r\varphi
_t(x)&=&\pmatrix{r\cr 2}
\frac{-\nabla_{1}^2\hat J(0)}{\hat I(0)} \hat\varphi(0)\sum
_{t=0}^\infty
m^t\sum_{x\in{\mathbb Z}^d}|x_1|^{r-2}\varphi_t(x)
\\[-2pt]
&&{}+O\biggl(\biggl(1-\frac{m}{m_{\rm c}}\biggr)^{-1-{r}/2+\epsilon}\biggr).
\end{eqnarray*}
Repeatedly using this recursion $j$ times and recalling $r-2j=q$, we
obtain
\begin{eqnarray*}
\sum_{t=0}^\infty m^t\sum_{x\in{\mathbb Z}^d}|x_1|^r\varphi
_t(x)&=&\prod_{i=0}^{j-1}
\pmatrix{r-2i\cr 2}\biggl(\frac{-\nabla_{1}^2\hat J(0)}{\hat I(0)} \hat
\varphi(0) \biggr)^j
\\[-2pt]
&&{}\times\sum_{t=0}^\infty m^t\sum_{x\in{\mathbb Z}^d
}|x_1|^{r-2j}\varphi_t(x)
+O\biggl(\biggl(1-\frac{m}{m_{\rm c}}\biggr)^{-1-{r}/2+\epsilon}\biggr)
\\[-2pt]
&=&\frac{\Gamma(r+1)}{2^j\Gamma(r-2j+1)}\biggl(\frac{-\nabla_{1}
^2\hat J(0)}
{\hat I(0)} \hat\varphi(0)\biggr)^j
\\[-2pt]
&&{}\times\sum_{t=0}^\infty m^t\sum_{x\in{\mathbb
Z}^d}|x_1|^q\varphi_t(x)
+O\biggl(\biggl(1-\frac{m}{m_{\rm c}}\biggr)^{-1-{r}/2+\epsilon}\biggr).
\end{eqnarray*}
Notice that, by (\ref{eq:chi-asy}) and (\ref{eq:ratio}),
\begin{eqnarray*}
\frac{-\nabla_{1}^2\hat J(0)}{\hat I(0)} \hat\varphi(0)=\frac
{2C_{\mathrm{II}}v_\alpha}
{1-{m/m_{\rm c}}}+O\biggl(\biggl(1-\frac{m}{m_{\rm c}}\biggr)^{-1+\epsilon
}\biggr)
\end{eqnarray*}
and that, by (\ref{eq:pr4r<2}) for $\alpha>2$ and (\ref{eq:Kr-evaluation}),
\begin{eqnarray*}
\sum_{t=0}^\infty m^t\sum_{x\in{\mathbb Z}^d}|x_1|^q\varphi
_t(x)=\Gamma
(q+1)
\frac{C_{\mathrm I}(C_{\mathrm{II}}v_\alpha)^{{q}/2}}{(1-
{m/m_{\rm
c}})^{1+{q}/2}}
+O\biggl(\biggl(1-\frac{m}{m_{\rm c}}\biggr)^{-1-{q}/2+\epsilon}\biggr).
\end{eqnarray*}
Therefore, we arrive at
\begin{eqnarray*}
&&\sum_{t=0}^\infty m^t\sum_{x\in{\mathbb Z}^d}|x_1|^r\varphi
_t(x)
\\
&&\qquad =\frac
{\Gamma(r+1)}{2^j
\Gamma(q+1)}\biggl(\frac{2C_{\mathrm{II}}v_\alpha}{1-{m/m_{\rm
c}}}\biggr)^j\Gamma(q+1)
\frac{C_{\mathrm I}(C_{\mathrm{II}}v_\alpha)^{{q}/2}}{(1-
{m/m_{\rm
c}})^{1+{q}/2}}
\\
&&{}\qquad\quad +O\biggl(\biggl(1-\frac{m}{m_{\rm c}}\biggr)^{-1-{r}/2+\epsilon}\biggr)
\\
&&\qquad =\Gamma(r+1) \frac{C_{\mathrm I}(C_{\mathrm{II}}v_\alpha)^{
{r}/2}}{(1-
{m/ m_{\rm c}})^{1+
{r}/2}}+O\biggl(\biggl(1-\frac{m}{m_{\rm c}}\biggr)^{-1-{r}/2+\epsilon
}\biggr).
\end{eqnarray*}
This completes the proof of Theorem~\ref{theorem:main1}.

\subsection{\texorpdfstring{Proof of
Theorem~\protect\ref{theorem:main2}}{Proof of Theorem~1.2}}

It is very easy to identify the main term for $\alpha\ne2$. First, by the
binomial expansion of the main term in (\ref{eq:main1}),
%
\begin{eqnarray}\label{eq:binom-exp}
\quad&&\biggl(1-\frac{m}{m_{\rm c}}\biggr)^{-1-{r/(\alpha\wedge
2)}}\nonumber
\\
&&\qquad =\sum_{t=0}^\infty
\frac{(-{r/(\alpha\wedge2)}-1)(-{r/(\alpha\wedge
2)}-2)\cdots
(-{r/(\alpha\wedge2)}-t)}{t!}\nonumber
\\
&&{}\hspace*{13pt}\qquad\quad\times \biggl(-\frac{m}{m_{\rm c}}
\biggr)^t
\\
&&\qquad =\sum_{t=0}^\infty\frac{\Gamma({r/(\alpha\wedge
2)}+t+1)}{t! \Gamma
({r/(\alpha\wedge2)}+1)}\biggl(\frac{m}{m_{\rm c}}
\biggr)^t\nonumber
\\
&&\qquad =\frac1{\Gamma({r/(\alpha\wedge2)}+1)}\sum_{t=0}^\infty
\biggl(\frac{m}
{m_{\rm c}}\biggr)^t\frac1{t!}\int_0^\infty x^{t+{r/(\alpha
\wedge2)}} e^{-x}\,
 \mathrm{d}x.\nonumber
\end{eqnarray}
Then, by the steepest descent method, we obtain that, for every $\beta
\in{\mathbb R}$,
\begin{eqnarray*}
\int_0^\infty x^{t+\beta} e^{-x} \,\mathrm{d}x\sim\sqrt{2\pi
(t+\beta)}\biggl(
\frac{t+\beta}{e}\biggr)^{t+\beta}\qquad\mbox{as }t\to\infty.
\end{eqnarray*}
Using this for $\beta=0,\frac{r}{\alpha\wedge2}$, we conclude that, as
$t\to\infty$,
\begin{eqnarray*}
\frac1{t!}\int_0^\infty x^{t+{r/(\alpha\wedge2)}} e^{-x}
\,\mathrm{d}x
&\sim&\biggl(\frac{t+{r/(\alpha\wedge2)}}t\biggr)^{t+
1/2}\biggl(\frac{t
+{r/(\alpha\wedge2)}}{e}\biggr)^{{r/(\alpha\wedge2)}}
\\[-2pt]
&\sim&
t^{{r/(\alpha\wedge2)}},
\end{eqnarray*}
which implies that the large-$t$ asymptotic expression for the
coefficient of
$m^t$ in~(\ref{eq:binom-exp}) is
$m_{\rm c}^{-t} t^{{r/(\alpha\wedge2)}}/\Gamma(\frac
{r}{\alpha\wedge2}+1)$,
hence the expression for the constant in (\ref{eq:main2}) for $\alpha
\ne2$.

There are many other ways to derive the above asymptotic expression.
One of
them is to notice that $x^te^{-x}/t!$ in (\ref{eq:binom-exp}) is the
probability
density for the sum of independent mean-one exponential random
variables. Then,
we use Jensen's inequality and apply the law of large numbers if
$\frac{r}{\alpha\wedge2}\le1,$ or exactly compute integer-power
moments for the
exponential random variables if $\frac{r}{\alpha\wedge2}>1$. We omit the
details.

To identify the main term for $\alpha=2$ in (\ref{eq:main2}), as well
as to
obtain the error estimates for all $\alpha>0$, we simply use
\cite{fo90}, Theorems~3A and 4. For convenience, we summarize a slightly
simplified version of these results as follows.

\begin{theorem}[(\cite{fo90}, Theorems~3A and 4)]\label{theorem:fo90}
\textup{(i)} Let
\begin{eqnarray*}
f(z)=(1-z)^{-1-\beta}\biggl(\log\frac1{1-z}\biggr)^\gamma,
\end{eqnarray*}
where $\beta\notin-\mathbb{N}\equiv{\mathbb Z}\setminus{\mathbb
Z}_+$ and $\gamma\notin
{\mathbb Z}_+$ are real or
complex numbers. Then, the coefficient $f_t$ of $f(z)=\sum_tf_tz^t$ satisfies
\begin{eqnarray*}
f_t\sim\frac{t^\beta(\log t)^\gamma}{\Gamma(1+\beta)}\qquad\mbox
{as }
t\to\infty.
\end{eqnarray*}\vspace*{-12pt}
\begin{longlist}[(ii)]
\item[(ii)]
Let $f(z)$ be analytic in $|z|<1$ and
\begin{eqnarray*}
f(z)=O(|1-z|^{-1-\beta})\qquad\mbox{as }z\to1
\end{eqnarray*}
for some real number $\beta>0$. Then, the coefficient $f_t$ of
$f(z)=\sum_tf_tz^t$ satisfies
\begin{eqnarray*}
f_t=O(t^\beta)\qquad\mbox{as }t\to\infty.
\end{eqnarray*}
\end{longlist}
\end{theorem}

The main term for $\alpha=2$ in (\ref{eq:main2}) is obtained by setting
$\beta=\gamma=\frac{r}2$ in Theorem~\ref{theorem:fo90}(i). For the error
estimates, we use Theorem~\ref{theorem:fo90}(ii) with $\beta=\frac
{r}2$ for
$\alpha=2$ and $\beta=\frac{r}{\alpha\wedge2}-\epsilon>0$ for
$\alpha\ne2$.
This completes the proof of Theorem~\ref{theorem:main2}.

\section{Proof of the key propositions}\label{s:key-prop}

In this section, we prove
Propositions~\ref{proposition:r<alpha&2} and \ref{proposition:2<r<alpha},
these being key propositions used in the previous section to prove the main
theorem. In Section~\ref{ss:prop2}, we first prove
Proposition~\ref{proposition:2<r<alpha}. Then, in Section~\ref{ss:prop1},
we use (\ref{eq:2<r<alpha-varphibd}) in Proposition~\ref
{proposition:2<r<alpha}
to show Proposition~\ref{proposition:r<alpha&2} for $\alpha>2$.

\subsection{\texorpdfstring{Proof of
Proposition~\protect\ref{proposition:2<r<alpha}}{Proof of Proposition~3.2}}
\label{ss:prop2}
Below, we prove Proposition~\ref{proposition:2<r<alpha} by using the results
already obtained in \cite{csI,csII,hhh,hhs08} and alternately applying
the following two
lemmas.

\begin{lemma}\label{lmm:induction1}
Let $\alpha>2$, $l\in\{1,2,\dots,\langle\alpha\rangle\}$ and
suppose that
(\ref{eq:2<r<alpha-IJbds}) holds for any $\nu\in\{0,1,\dots,l\vee
2\}$ and
(\ref{eq:2<r<alpha-varphibd}) holds for any $n\in\{0,\dots,l-1\}$.
Then, (\ref{eq:2<r<alpha-varphibd}) holds for $n=l$.
\end{lemma}

\begin{lemma}\label{lmm:induction2}
Let $\alpha>2$ and suppose that (\ref{eq:2<r<alpha-varphibd}) holds
for $n=2l$,
where $l\in\{1,\dots,\langle\frac\alpha2\rangle\}$ (note that
$\alpha-2\le2\langle\frac\alpha2\rangle<\alpha$). Then,
(\ref{eq:2<r<alpha-IJbds}) holds for any $\nu\in(n,n+2]$ if
$n+2<\alpha$,
or for any $\nu\in(n,\alpha)$ if $\alpha\le n+2$.
\end{lemma}

We will prove these lemmas after completing the proof of
Proposition~\ref{proposition:2<r<alpha}. For random walk,
(\ref{eq:2<r<alpha-IJbds}) always holds as mentioned earlier and we
therefore only need Lemma~\ref{lmm:induction1}.

We now begin by proving Proposition~\ref{proposition:2<r<alpha}.
First, we note that (\ref{eq:2<r<alpha-IJbds}) for $\nu\in[0,2]$ and
(\ref{eq:2<r<alpha-varphibd}) for $n=0$ have been proven in the current
setting \mbox{\cite{csI,csII,hhh,hhs08}}; the result in \cite{hhs08} for
self-avoiding walk is only valid at $\theta=0$. However, it is not
hard to extend
the result to nonzero $\theta$ by splitting the denominator in
(\ref{eq:lace-exp-solution}) into $1-\hat J(k,m)$ and
$\hat J(k,m)-\hat J(k,me^{i\theta}),$ and estimating the latter as
$m^t-(me^{i\theta})^t=m^t(1-e^{i\theta})\sum_{s=0}^{t-1}e^{i\theta s}$
[which equals $O(\theta)tm^t$ for $|\theta|\ll1$]. We omit the
details. Then, by
Lemma~\ref{lmm:induction1} with $l=1$, we obtain (\ref{eq:2<r<alpha-varphibd})
for $n=1$. With this conclusion and again using Lemma~\ref{lmm:induction1},
but now with $l=2$, we obtain (\ref{eq:2<r<alpha-varphibd}) for $n=2$.
With this
conclusion and using Lemma~\ref{lmm:induction2}, we further obtain
(\ref{eq:2<r<alpha-IJbds}) for $\nu\in(2,4]$ or $\nu\in(2,\alpha)$,
depending on
whether $\alpha>4$ or $\alpha\le4$. We can repeat this, using
Lemmas~\ref{lmm:induction1} and \ref{lmm:induction2} alternately,
until $n$
reaches $\langle\alpha\rangle$. Let $\tilde l=\langle\frac\alpha
2\rangle$.
We see that
\begin{eqnarray*}
&&
\left.\begin{array}{l}
(\ref{eq:2<r<alpha-IJbds})_{\nu\in[0,2]}\\[5pt]
(\ref{eq:2<r<alpha-varphibd})_{n=0}
\end{array}
\right\}
\stackrel{\mathrm{Lemma~\scriptsize{\ref{lmm:induction1}}}}\Longrightarrow
(\ref{eq:2<r<alpha-varphibd})_{n=1,2}
\stackrel{\mathrm{Lemma~\scriptsize{\ref{lmm:induction2}}}}\Longrightarrow
(\ref{eq:2<r<alpha-IJbds})_{\nu\in(2,4]}
\stackrel{\mathrm{Lemma~\scriptsize{\ref{lmm:induction1}}}}\Longrightarrow
\cdots
\\
&&\hspace*{67pt}\stackrel{\mathrm{Lemma~\scriptsize{\ref
{lmm:induction1}}}}\Longrightarrow
(\ref{eq:2<r<alpha-varphibd})_{n=2\tilde l-1,2\tilde l}
\\
&&\hspace*{67pt}\stackrel{\mathrm{Lemma~\scriptsize{\ref{lmm:induction2}}}}\Longrightarrow
(\ref{eq:2<r<alpha-IJbds})_{\nu\in(2\tilde l,\alpha)}
\bigl(\mathop{\Longrightarrow}\limits^{\mathrm{Lemma~\scriptsize{\ref{lmm:induction1}}}}_{\mathrm{if }\ \alpha>2\tilde l+1}
(\ref{eq:2<r<alpha-varphibd})_{n=2\tilde l+1}\bigr).
\end{eqnarray*}
This completes the proof of Proposition~\ref{proposition:2<r<alpha}.

\begin{pf*}{Proof of Lemma~\ref{lmm:induction1}}
First, by using (\ref{eq:2<r<alpha-IJbds}) for $\nu=2$ and
(\ref{eq:2<r<alpha-varphibd}) for $n=0$,
we prove $|\nabla_{1}\hat I(k)|\le O(|\hat\varphi(k)|^{-1/2})$; the
proof of
$|\nabla_{1}\hat J(k)|\le O(|\hat\varphi(k)|^{-1/2})$ is almost
identical and
thus we omit it. By the ${\mathbb Z}^d$-symmetry of the models
and using\vadjust{\goodbreak} $|\sin(k_1x_1)|\le|k_1x_1|$ and (\ref{eq:2<r<alpha-IJbds})
for $\nu=2$,
we obtain
\begin{eqnarray*}
|\nabla_{1}\hat I(k)|&=&\Biggl|\sum_{t=0}^\infty m^t\sum_{x\in
{\mathbb Z}^d
}x_1\sin(k_1x_1)
 I_t(x) e^{i(k_2x_2+\cdots+k_dx_d)}\Biggr|
 \\[-2pt]
&\le&|k_1|\sum_{t=0}^\infty m^t\sum_{x\in{\mathbb Z}^d}x_1^2
|I_t(x)|\le O(|k_1|).
\end{eqnarray*}
However, by (\ref{eq:2<r<alpha-varphibd}) for $n=0$, we have
$|\hat\varphi(k)|\le O(|k|^{-2})$, which implies that
$|k_1|\le|k|\le O(|\hat\varphi(k)|^{-1/2})$, as required.

We now use this bound to complete the proof of Lemma~\ref{lmm:induction1}.
First, by
differentiating (\ref{eq:lace-exp-fourier}) and solving the resulting equation
for $\nabla_{1}^l\hat\varphi(k)$, we have that, for $l\in\mathbb{N}$,\vspace*{-2pt}
\begin{eqnarray*}
\nabla_{1}^l\hat\varphi(k)&=&\nabla_{1}^l\hat I(k)+\sum
_{j=0}^l\pmatrix
{l\cr j}\nabla_{1}^j
\hat J(k) \nabla_{1}^{l-j}\hat\varphi(k)
\\[-2pt]
&=&\frac1{1-\hat J(k)}\Biggl(\nabla_{1}^l\hat I(k)+\sum_{j=1}^l\pmatrix
{l\cr j}\nabla_{1}^j
\hat J(k) \nabla_{1}^{l-j}\hat\varphi(k)\Biggr).\vspace*{-2pt}
\end{eqnarray*}
By (\ref{eq:lace-exp-solution}), (\ref{eq:Ihat}) and (\ref
{eq:property1}), we have
$|1-\hat J(k)|^{-1}=O(|\hat\varphi(k)|)$. By (\ref
{eq:2<r<alpha-IJbds}) for
$\nu\ge2$ or using $|\nabla_{1}\hat I(k)|\le O(|\hat\varphi(k)|^{-1/2})$,
we obtain\vspace*{-2pt}
\begin{eqnarray*}
\biggl|\frac{\nabla_{1}^l\hat I(k)}{1-\hat J(k)}\biggr|\le O
(|\hat
\varphi(k)|
)\times
\cases{
|\hat\varphi(k)|^{-1/2},&\quad$l=1$,
\cr
\biggl(1-\dfrac{m}{m_{\rm c}}\biggr)^{1-{l/2}+\epsilon},&\quad$l=2,\dots
,\langle\alpha
\rangle$,
}\vspace*{-2pt}
\end{eqnarray*}
which, by (\ref{eq:2<r<alpha-varphibd}) for $n=0$, is smaller than the
bound in
(\ref{eq:2<r<alpha-varphibd}) for $n=l$, yielding an error term. For $j=1,2$,
we also use (\ref{eq:2<r<alpha-varphibd}) for $n\le l-1$ to obtain\vspace*{-2pt}
\begin{eqnarray*}
&&\biggl|\frac{\nabla_{1}^j\hat J(k)}{1-\hat J(k)} \nabla
_{1}^{l-j}\hat
\varphi(k)
\biggr|
\\[-2pt]
&&\qquad \le O(|\hat\varphi(k)|^{j/2})
\\[-2pt]
&&{}\qquad\quad\times
\cases{
\biggl(1-\dfrac{m}{m_{\rm c}}+|\theta|+|k|^2\biggr)^{-1-{(l-j)}/2},
&\quad$l=j,j+1$,
\cr
\displaystyle\dfrac{(1-{m/m_{\rm c}})^{1-
{(l-j)}/2}}{(1-
{m/m_{\rm c}}+|\theta|
+|k|^2)^2},&\quad$l=j+2,\dots,\langle\alpha\rangle$,
}\vspace*{-2pt}
\end{eqnarray*}
which, again by (\ref{eq:2<r<alpha-varphibd}) for $n=0$, obeys the
required bound
in (\ref{eq:2<r<alpha-varphibd}) for $n=l$. Finally, for $j\ge3$ (hence
for $l\ge3$),\vspace*{-2pt}
\begin{eqnarray*}
&&\biggl|\frac{\nabla_{1}^j\hat J(k)}{1-\hat J(k)} \nabla
_{1}^{l-j}\hat
\varphi(k)
\biggr|
\\[-2pt]
&&\qquad\le\frac{O((1-{m/m_{\rm c}})^{1-{j}/2+\epsilon
})}{1-{m/m_{\rm c}}
+|\theta|+|k|^2}
\\
&&\qquad\quad{}\times
\cases{
\biggl(1-\dfrac{m}{m_{\rm
c}}+|\theta|+|k|^2\biggr)^{-1-{(l-j)}/2}&\quad$[l=j,j+1]$,
\cr
\displaystyle\dfrac{(1-{m/m_{\rm c}})^{1-
{(l-j)}/2}}{(1-
{m/m_{\rm c}}+|\theta|
+|k|^2)^2}&\quad$[l=j+2,\dots,\langle\alpha\rangle]$
}
\\
&&\qquad\le\frac{O((1-{m/m_{\rm c}})^{1-{l}/2+\epsilon
})}{(1-{m/m_{\rm c}}
+|\theta|+|k|^2)^2},
\end{eqnarray*}
which is smaller [by the factor $(1-\frac{m}{m_{\rm c}})^\epsilon$]
than the bound
in (\ref{eq:2<r<alpha-varphibd}), yielding\vspace*{1pt} an error term. This
completes the
proof of Lemma~\ref{lmm:induction1}.
\end{pf*}

\begin{pf*}{Proof of Lemma~\ref{lmm:induction2}}
First, we recall (\ref{eq:I-def}) and (\ref{eq:J-def}). Since\break
$\sum_x|x_1|^\nu D(x)<\infty$ provided that $\nu<\alpha$, (\ref
{eq:2<r<alpha-IJbds})
always holds for random walk. Moreover, for oriented percolation, there
is a
constant $C_\nu<\infty$ such that
\begin{eqnarray*}
\sum_{t=0}^\infty m^t\sum_{x\in{\mathbb Z}^d}|x_1|^\nu
|J_t^\mathrm{OP}(x)|&\le&
p\sum_{t=1}^\infty
m^t\sum_{x,y\in{\mathbb Z}^d}|y_1+x_1-y_1|^\nu|\pi
_{t-1}^\mathrm{OP}(y)| D(x-y)
\\
&\le& C_\nu pm\sum_{t=1}^\infty m^{t-1}\sum_{y\in{\mathbb
Z}^d}(|y_1|^\nu
+1)|\pi_{t-1}^\mathrm{OP}
(y)|,
\end{eqnarray*}
where we have used the fact that, for any $a_1,\dots,a_n\in{\mathbb R}$,
%
\begin{eqnarray}\label{eq:naive}
\Biggl|\sum_{j=1}^na_j\Biggr|^\nu\le\Bigl(n\max_{1\le j\le
n}|a_j|\Bigr)^\nu
=n^\nu\max_{1\le j\le n}|a_j|^\nu\le n^\nu\sum_{j=1}^n|a_j|^\nu.
\end{eqnarray}
Since $\sum_{s=0}^\infty m^s\sum_{y\in{\mathbb Z}^d}|\pi
_s^\mathrm{OP}(y)|=O(1)$
uniformly in
$m\le m_{\rm c}$ \cite{csI}, it suffices to show that, for
self-avoiding walk and oriented percolation, (\ref
{eq:2<r<alpha-varphibd}) for
$n=2l$, where $l\in\{1,\dots,\langle\frac\alpha2\rangle\}$ implies
that
%
\begin{eqnarray}\label{eq:2<r<alpha-pibd}
\sum_{t=0}^\infty m^t\sum_{x\in{\mathbb Z}^d}|x_1|^{2l+q}|\pi
_t(x)|\le
O\biggl(\biggl(1-\frac{m}{m_{\rm c}}\biggr)^{1-{(2l+q)}/2+\epsilon}\biggr)
\end{eqnarray}
for any $q\in(0,2]$ if $2l+2<\alpha$, or for any $q\in(0,\alpha
-2l)$ if
$\alpha\le2l+2$.

As we mentioned earlier, $\pi_t(x)$ is an alternating sum of the
lace expansion coefficients. More precisely,
\begin{eqnarray*}
\pi_t(x)=\sum_{N=0}^\infty(-1)^N\pi_t^{(N)}(x),
\end{eqnarray*}
where $\pi_t^{(N)}(x)\ge0$ is the model-dependent $N$th
expansion coefficient (see, e.g., \cite{csI,s06} for the precise
definitions of the expansion coefficients). Due to the subadditivity argument
for self-avoiding walk and by the BK inequality \cite{bk}
for percolation,
it is known that the expansion coefficients satisfy the following diagrammatic
bounds, in which each line corresponds to a 2-point function. For
self-avoiding walk,
\begin{eqnarray}\label{eq:saw-diagbds}
\pi_t^{(0)}(x)&\equiv&0 ,\qquad
\pi_t^{(1)}(x)\le {}\mbox{
\includegraphics{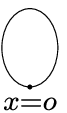}
},\nonumber
\\[-8pt]\\[-8pt]
\pi_t^{(2)}(x)&\le& {}\mbox{
\includegraphics{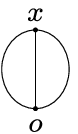}
},\qquad
\pi_t^{(3)}(x)\le {}\mbox{
\includegraphics{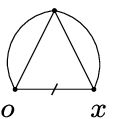}
},\nonumber
\end{eqnarray}
where the bounding diagram for $\pi_t^{(1)}(x)$ is the $t$-step
self-avoiding loop at $x=o$, hence proportional to $\delta_{x,o}$, and the
diagram for $\pi_t^{(2)}(x)$ is the product of three 2-point functions
$\varphi_s^\mathrm{SAW}(x) \varphi_{s'}^\mathrm
{SAW}(x) \varphi_{s''}^\mathrm{SAW}(x)$
summed over
all possible combinations of $s,s',s''\in\mathbb{N}$ satisfying $s+s'+s''=t$,
and so
on. The unlabeled vertices in the diagrams for $\pi_t^{(3)}(x)$ and
the higher order expansion coefficients are summed over ${\mathbb
Z}^d$. For oriented
percolation,
%
\begin{eqnarray}\label{eq:op-diagbds}
\pi_t^{(0)}(x)\le{}\mbox{
\includegraphics{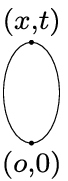}
},\qquad
\pi_t^{(1)}(x)\le{}\mbox{
\includegraphics{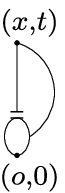}
},\qquad
\pi_t^{(2)}(x)\le{}\mbox{
\includegraphics{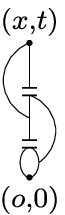}
} +
{}\mbox{
\includegraphics{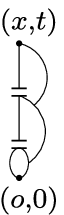}
},
\end{eqnarray}
where the bounding diagram for $\pi_t^{(0)}(x)$ is $\varphi
_t^\mathrm{OP}(x)^2$
and that for $\pi_t^{(1)}(x)$ is the product of five 2-point functions
concatenated in the depicted way, and so on. The upward direction of the
diagrams is the time-increasing direction and the unlabeled vertices are
summed over space--time ${\mathbb Z}^d\times{\mathbb Z}_+$. For more
details, we refer to
\cite{s07}.

First, we prove (\ref{eq:2<r<alpha-pibd}) for self-avoiding walk. Since
$\pi_t^{(0)}(x)\equiv0$ and $\pi_t^{(1)}(x)\propto\delta_{x,o}$,
it suffices to investigate the contributions from $\pi_t^{
(N)}(x)$ for
$N\ge2$. For $\pi_t^{(2)}(x)$, since
%
\begin{eqnarray}\label{eq:pi2sawbd}
\pi_t^{(2)}(x)\le\mathop{\sum_{s,s',s''\in\mathbb{N}}}_{(s+s'+s''=t)}
\varphi_s^\mathrm{SAW}(x) \varphi_{s'}^\mathrm
{SAW}(x) \varphi_{s''}^\mathrm{SAW}(x),
\end{eqnarray}
we obtain
\begin{eqnarray*}
\sum_{t=0}^\infty m^t\sum_{x\in{\mathbb Z}^d}|x_1|^{2l+q} \pi_t^{
(2)}(x)&\le&
\biggl(\sum_{x\in{\mathbb Z}^d}|x_1|^q\sum_{s,s'\in\mathbb
{N}}m^{s+s'}\varphi
_s^\mathrm{SAW}(x)
\varphi_{s'}^\mathrm{SAW}(x)\biggr)
\\
&&{}\times\biggl(\sup_{x\in{\mathbb Z}^d}|x_1|^{2l}\sum_{s''\in
\mathbb{N}}m^{s''}
\varphi_{s''}^\mathrm{SAW}(x)\biggr)
\\
&\le& B^{(q)}W^{(2l)},
\end{eqnarray*}
where
\begin{eqnarray*}
B^{(\nu)}&=&\sup_{y\in{\mathbb Z}^d}\sum_{x\in{\mathbb
Z}^d}|x_1|^\nu\sum
_{t\in\mathbb{N}}m^t
\varphi_t^\mathrm{SAW}(x)\sum_{s=0}^\infty m^s\varphi
_s^\mathrm{SAW}(y-x),
\\
W^{(\nu)}&=&\sup_{x\in{\mathbb Z}^d}|x_1|^\nu\sum_{t\in\mathbb
{N}}m^t\varphi
_t^\mathrm{SAW}(x).
\end{eqnarray*}
Similarly to the above and the derivation of \cite{hhh}, formula (2.42), by using
(\ref{eq:naive}) and diagrammatic bounds of the form (\ref{eq:saw-diagbds}),
we can show that
\begin{eqnarray}\label{eq:piNsawbd}
&& \sum_{t=0}^\infty m^t\sum_{x\in{\mathbb Z}^d}|x_1|^{2l+q} \pi_t^{
(N)}(x)\nonumber
\\[-8pt]\\[-8pt]
&&\qquad\le
N^{2l+q+2}\bigl(B^{(0)}\bigr)^{N-2}B^{(q)}W^{(2l)}\qquad[N\ge2].\nonumber
\end{eqnarray}
It is immediate from the definition (\ref{eq:2pt-def}) that
$\varphi_t^\mathrm{SAW}(x)\le\delta_{x,o}\delta
_{t,0}+(D*\varphi
_{t-1}^\mathrm{SAW})(x)$.
By this, we have
%
\begin{eqnarray}\label{eq:B0sawbd}
B^{(0)}&\le&\sup_{y\in{\mathbb Z}^d}\sum_{x\in{\mathbb Z}^d}\sum
_{t\in\mathbb{N}
}m^t\varphi_t^\mathrm{SAW}(x)
\biggl(\delta_{x,y}+\sum_{s\in\mathbb{N}}m^s(D*\varphi
_{s-1}^\mathrm{SAW}
)(y-x)\biggr)\nonumber
\\
&\le& W^{(0)}+\sup_{y\in{\mathbb Z}^d}\sum_{x\in{\mathbb Z}^d}\sum
_{t\in\mathbb{N}
}m^t(D*\varphi_{t-
1}^\mathrm{SAW})(x)\sum_{s\in\mathbb{N}}m^s(D*\varphi
_{s-1}^\mathrm{SAW})(y-x)
\\
&\le& W^{(0)}+m^2\int_{[-\pi,\pi]^d}\frac{\mathrm{d}^dk}{(2\pi
)^d} \hat
D(k)^2|\hat\varphi^\mathrm{SAW}(k,m)|^2\nonumber
\end{eqnarray}
and
%
\begin{eqnarray}
W^{(0)}&\le&\sup_{x\in{\mathbb Z}^d}\sum_{t\in\mathbb
{N}}m^t(D*\varphi
_{t-1}^\mathrm{SAW})(x)\nonumber
\\[-2pt]
&\le& m\|D\|_\infty+\sup_{x\in{\mathbb Z}^d}\sum_{t=2}^\infty
m^t(D*D*\varphi
_{t-2}^\mathrm{SAW})
(x)
\\[-2pt]
&\le& m\|D\|_\infty+m^2\int_{[-\pi,\pi]^d}\frac{\mathrm{d}^dk}{(2\pi
)^d} \hat
D(k)^2|\hat\varphi^\mathrm{SAW}(k,m)|.\nonumber
\end{eqnarray}
By (\ref{eq:2<r<alpha-varphibd}) for $n=0$ and $\|D\|_\infty
=O(L^{-d})$, we
can show that $B^{(0)}=O(L^{-d})$ uniformly in $m\le m_{\rm c}$ if $d>4$,
hence the summability of (\ref{eq:piNsawbd}) over $N\ge2$ when $L\gg1$.
Moreover,
by (\ref{eq:2<r<alpha-varphibd}) for $n=2l$,
%
\begin{eqnarray}\label{eq:W2lsawbd}
W^{(2l)}&\le&\int_{[-\pi,\pi]^d}\frac{\mathrm{d}^dk}{(2\pi)^d}
|\nabla_{1}^{2l}
\hat\varphi^\mathrm{SAW}(k,m)|\nonumber
\\[-2pt]
&\le&
O\biggl(\biggl(1-\frac
{m}{m_{\rm c}}\biggr)^{1-l}
\biggr)\int_{[-\pi,\pi]^d}\frac{\mathrm{d}^dk}{|k|^4}
\\[-2pt]
&\stackrel{d>4}=& O\biggl(\biggl(1-\frac{m}{m_{\rm
c}}\biggr)^{1-l}\biggr).\nonumber
\end{eqnarray}
Therefore,
\begin{eqnarray*}
\sum_{N=2}^\infty\sum_{t=0}^\infty m^t\sum_{x\in{\mathbb
Z}^d}|x_1|^{2l+q}
\pi_t^{
(N)}(x)\le O\biggl(\biggl(1-\frac{m}{m_{\rm c}}\biggr)^{1-l}\biggr)B^{(q)}.
\end{eqnarray*}

To complete the proof of (\ref{eq:2<r<alpha-pibd}) for self-avoiding
walk, it
suffices to show that there is an $\epsilon>0$ such that
$B^{(q)}=O((1-\frac{m}{m_{\rm c}})^{-{q}/2+\epsilon})$. For $q=2$,
we use (\ref{eq:2<r<alpha-varphibd}) for $n=2$ and take an arbitrary
$\epsilon\in(0,1\wedge\frac{d-4}2)$ to obtain
\begin{eqnarray*}
B^{(2)}&\le&\int_{[-\pi,\pi]^d}\frac{\mathrm{d}^dk}{(2\pi)^d}
|\hat
\varphi^\mathrm{SAW}(k,m) \nabla_{1}^2\hat\varphi
^\mathrm{SAW}(k,m)|
\\
&\le& O\biggl(\biggl(1-\frac{m}{m_{\rm c}}\biggr)^{-1+\epsilon}\biggr)\int_{[-\pi
,\pi]^d}
\frac{\mathrm{d}^dk}{|k|^{2(2+\epsilon)}}\le O\biggl(\biggl(1-\frac
{m}{m_{\rm c}}\biggr)^{-1
+\epsilon}\biggr).
\end{eqnarray*}
For $q\in(0,2)$, we first note that
%
\begin{eqnarray}\label{eq:Bqprebd}
\quad B^{(q)}\le\frac1{K_q}\int_0^\infty\frac{\mathrm
{d}u}{u^{1+q}}\int_{[-\pi,
\pi]^d}\frac{\mathrm{d}^dk}{(2\pi)^d} |\hat\varphi^\mathrm
{SAW}(k,m)
\bar\Delta_{\vec
u}\hat\varphi^\mathrm{SAW}(k,m)|.
\end{eqnarray}
It is known that, by \cite{hhs08}, Proposition~2.6, with an improvement
due to the same argument as in \cite{csII}, Proposition~2.1,
\begin{eqnarray*}
&&|\bar\Delta_{\vec u}\hat\varphi^\mathrm{SAW}(k,m)|
\\
&&\qquad \le\sum_{(j,j')=(0,\pm1),(1,-1)}\frac{O(1-\hat D(\vec
u))
}{1-{m/m_{\rm c}}
+1-\hat D(k+j\vec u)}
\\
&&\hspace*{98pt}{}\quad{}\times \frac{1}{1-{m/m_{\rm c}}+1-\hat D(k+j'\vec
u)}
\end{eqnarray*}
holds in the current setting, where the $O(1-\hat D(\vec u))$ term is
uniform in $k\in[-\pi,\pi]^d$ and $m\le m_{\rm c}$. Substituting
this, and
(\ref{eq:2<r<alpha-varphibd}) for $n=0,$ into (\ref{eq:Bqprebd}), and
using the
translation invariance and the ${\mathbb Z}^d$-symmetry of $D$ and the Schwarz
inequality (see \cite{csII}, formulas (4.27)--(4.29)), we end up with
\begin{eqnarray*}
B^{(q)}
&\le&\int_0^\infty\frac{\mathrm{d}u}{u^{1+q}}\int_{[-\pi,\pi]^d}\frac{\mathrm{d}^dk}{(2\pi)^d}
  \frac{O(1-\hat D(\vec u))}{(1-
{m/m_{\rm c}}
+1-\hat D(k))^2}
\\
&&{}\hspace*{81pt}\times\frac{1}{1-{m}/{m_{\rm c}}+1-\hat D(k-\vec
u)}
\\
&\le&\int_0^\infty\mathrm{d}u \frac{1-\hat D(\vec u)}{u^{1+q}}\int
_{[-\pi,
\pi]^d}\frac{\mathrm{d}^dk}{(2\pi)^d} \frac{O((1-{m/m_{\rm
c}})^{-{q}/2
+\epsilon})}{(1-\hat D(k))^{2-{q}/2+\epsilon}(1-\hat D(k-\vec u))}
\end{eqnarray*}
for any $\epsilon\in(0,\frac{q}2)$. However, by following the proof of
\cite{csII}, formula (4.30), we can show that
\begin{eqnarray*}
\int_{[-\pi,\pi]^d}\frac{\mathrm{d}^dk}{(2\pi)^d} \frac1{(1-\hat D(k))^{2
-{q}/2+\epsilon}(1-\hat D(k-\vec u))}\le O
\bigl(u^{(d-6+q-2\epsilon)
\wedge0}\bigr),
\end{eqnarray*}
hence
\begin{eqnarray}\label{eq:Bqbd}
\quad B^{(q)}&\le& O\biggl(\biggl(1-\frac{m}{m_{\rm c}}\biggr)^{-{q}/2+\epsilon
}\biggr)\biggl(
\int_0^1\frac{\mathrm{d}u}u u^{(d-4-2\epsilon)\wedge(2-q)}+\int
_1^\infty
\frac{\mathrm{d}u}{u^{1+q}}\biggr)\nonumber
\\[-9pt]\\[-9pt]
&=&O\biggl(\biggl(1-\frac{m}{m_{\rm c}}\biggr)^{-{q}/2+\epsilon}
\biggr)\nonumber
\end{eqnarray}
if $\epsilon<\frac{d-4}2$. This completes the proof of (\ref
{eq:2<r<alpha-pibd})
for self-avoiding walk.

For oriented percolation, similarly to the proof of \cite{csII},
Lemma~3, by
using (\ref{eq:naive}) and diagrammatic bounds of the form (\ref
{eq:op-diagbds}),
we can show that, for $N\ge0$,
%
\begin{eqnarray}\label{eq:piNopbd}
\qquad &&\sum_{t=0}^\infty m^t\sum_{x\in{\mathbb Z}^d}|x_1|^{2l+q} \pi_t^{
(N)}(x)\nonumber
\\[-2pt]
&&\qquad\le(N+1)^{2l+q}\bigl(T^{(0)}\bigr)^{N-2}\bigl(\bigl(N\bigl(1+T^{(0)}\bigr)
+T^{(0)}\bigr) T^{(0)} V^{(q)}
\\[-2pt]
&&\hspace*{102pt}\qquad\quad{}+N\bigl((N-1)\bigl(1+T^{(0)}\bigr)+3T^{(0)}\bigr)T^{
(q)}
V^{(0)}\bigr),\nonumber
\end{eqnarray}
where
\begin{eqnarray*}
V^{(\nu)}
&=&\sup_{(x,t)\in{\mathbb Z}^{d+1}}\sum_{(y,s)\in{\mathbb Z}^{d+1}}|y_1|^{2l}
(mD*\varphi_s^\mathrm{OP})(y) m^s |y_1-x_1|^\nu
\\[-2pt]
&&{}\hspace*{81pt}\times (D*\varphi
_{s-t}^\mathrm{OP}
)(y-x),
\\[-2pt]
T^{(\nu)}
&=&\sup_{(x,t)\in{\mathbb Z}^{d+1}}\sum_{(y,s),(y',s')\in{\mathbb Z}^{d+1}}
(mD*\varphi_s^\mathrm{OP})(y) m^s |y_1-x_1|^\nu
\\
&&{}\hspace*{105pt}\times(D*\varphi
_{s'-t}^\mathrm{OP}
)(y'-x)
\\[-2pt]
&&{}\hspace*{105pt}\times\bigl(\varphi_{s-s'}^\mathrm
{OP}(y-y')+\varphi_{s'-s}^\mathrm{OP}
(y'-y)\bigr).
\end{eqnarray*}
Notice that
%
\begin{eqnarray}\label{eq:T0bd}
\qquad\quad T^{(0)}\le2m\int_{[-\pi,\pi]^d}\frac{\mathrm{d}^dk}{(2\pi
)^d} \hat D(k)^2
\int_{-\pi}^\pi\frac{\mathrm{d}\theta}{2\pi} |\hat\varphi
^\mathrm{OP}
(k,me^{i\theta})|
|\hat\varphi^\mathrm{OP}(k,e^{i\theta})|^2.
\end{eqnarray}
Using (\ref{eq:2<r<alpha-varphibd}) for $n=0$ and $\|D\|_\infty
=O(L^{-d})$, we
can show that $T^{(0)}=O(L^{-d})$ uniformly in $m\le m_{\rm c}$ if
$d>4$ and
$L\gg1$, hence the summability of (\ref{eq:piNopbd}) over $N\ge0$. Moreover,
by (\ref{eq:2<r<alpha-varphibd}) for $n=0,2l$ and using $|\hat
D(k)|\le1$,
we have
%
\begin{eqnarray}\label{eq:V0bd}
V^{(0)}
&\le&2^{2l}m\biggl(\int_{[-\pi,\pi]^d}
\frac{\mathrm{d}^dk}{(2\pi)^d}\int_{-\pi}^\pi\frac{\mathrm{d}\theta
}{2\pi}
|\nabla_{1}^{2l}\hat\varphi^\mathrm{OP}(k,me^{i\theta
})||\hat\varphi^\mathrm{OP}(k,e^{i
\theta})|\nonumber
\\[-2pt]
&&{}\hspace*{26pt}+\sum_{x\in{\mathbb Z}^d}|x_1|^{2l}D(x)\int_{[-\pi,\pi
]^d}\frac
{\mathrm{d}^dk}
{(2\pi)^d}\int_{-\pi}^\pi\frac{\mathrm{d}\theta}{2\pi} |\hat
\varphi^\mathrm{OP}
(k,me^{i\theta})|
\\[-2pt]
&&\hspace*{182pt}{}\hspace*{26pt}\times |\hat\varphi^\mathrm{OP}(k,e^{i\theta
})|\biggr)\nonumber
\\[-2pt]
&\stackrel{d>4}=&O\biggl(\biggl(1-\frac{m}{m_{\rm c}}\biggr)^{1-l}\biggr).\nonumber
\end{eqnarray}
To complete the proof of (\ref{eq:2<r<alpha-pibd}), it thus suffices to
show that
there is an $\epsilon>0$ such that
\begin{eqnarray*}
T^{(q)}=O\biggl(\biggl(1-\frac{m}{m_{\rm c}}\biggr)^{-{q}/2+\epsilon
}\biggr),\qquad
V^{(q)}=O\biggl(\biggl(1-\frac{m}{m_{\rm c}}\biggr)^{1-l-{q}/2+\epsilon
}\biggr).
\end{eqnarray*}

Here, we only explain the proof of the bound on $V^{(2)}$; the
bound on
$T^{(2)}$ can be proven quite similarly and the bounds on $T^{(q)}$
and $V^{(q)}$ for $q\in(0,2)$ can be proven by following a
similar line
of argument from (\ref{eq:Bqprebd}) through to (\ref{eq:Bqbd}). To prove
the bound
on $V^{(2)}$, we first note that
\begin{eqnarray}\label{eq:V2prebd}
\quad V^{(2)}&\le&2^{2l+2}m\biggl(\int_{[-\pi,\pi]^d}\frac{\mathrm
{d}^dk}{(2\pi)^d}
\int_{-\pi}^\pi\frac{\mathrm{d}\theta}{2\pi} |\nabla
_{1}^{2l}\hat
\varphi^\mathrm{OP}(k,
me^{i\theta})||\nabla_{1}^2\hat\varphi^\mathrm
{OP}(k,e^{i\theta})|\nonumber
\\
&&{}\hspace*{37pt}+\sigma^2\int_{[-\pi,\pi]^d}\frac{\mathrm{d}^dk}{(2\pi
)^d}\int_{-
\pi}^\pi\frac{\mathrm{d}\theta}{2\pi} |\nabla_{1}^{2l}\hat
\varphi
^\mathrm{OP}(k,me^{i
\theta})||\hat\varphi^\mathrm{OP}(k,e^{i\theta})|\nonumber
\\
&&{}\hspace*{37pt}+\sum_{x\in{\mathbb Z}^d}x_1^{2l}D(x)\int_{[-\pi,\pi
]^d}\frac{\mathrm{d}^dk}
{(2\pi)^d}\int_{-\pi}^\pi\frac{\mathrm{d}\theta}{2\pi} |\hat
\varphi^\mathrm{OP}(k,
me^{i\theta})|\nonumber
\\[-8pt]\\[-8pt]
&&{}\hspace*{197pt}\times |\nabla_{1}^2\hat\varphi^\mathrm
{OP}(k,e^{i\theta})|\nonumber
\\
&&{}\hspace*{37pt}+\sigma^2\sum_{x\in{\mathbb Z}^d}x_1^{2l}D(x)\int_{[-\pi
,\pi]^d}
\frac{\mathrm{d}^dk}{(2\pi)^d}\int_{-\pi}^\pi\frac{\mathrm{d}\theta
}{2\pi}
|\hat\varphi^\mathrm{OP}(k,me^{i\theta})|\nonumber
\\
&&{}\hspace*{37pt}\hspace*{186pt}\times |\hat\varphi
^\mathrm{OP}(k,e^{i\theta
})|\biggr).\nonumber
\end{eqnarray}
It is immediate from (\ref{eq:2<r<alpha-varphibd}) for $n=0,2$ that the
last two
lines are both $O(1)$ for $d>4$. Moreover, by (\ref{eq:V0bd}), the
second line
is $O((1-\frac{m}{m_{\rm c}})^{1-l})$ for $d>4$. For the first line,
we use the
following bounds due to (\ref{eq:2<r<alpha-varphibd}) for $n=2,2l$:
for any
$\epsilon\in(0,1)$,
\begin{eqnarray*}
|\nabla_{1}^{2l}\hat\varphi^\mathrm{OP}(k,me^{i\theta})|\le
\frac{O
((1-{m}/
{m_{\rm c}})^{-l+\epsilon})}{(|\theta|+|k|^2)^{1+\epsilon}},\qquad
|\nabla_{1}^2\hat\varphi^\mathrm{OP}(k,e^{i\theta})|\le O(|k|^{-4}),
\end{eqnarray*}
where the $O((1-\frac{m}{m_{\rm c}})^{-l+\epsilon})$ term is uniform in
$(k,\theta)\in[-\pi,\pi]^{d+1}$ and the $O(|k|^{-4})$ term is
uniform in
$\theta\in[-\pi,\pi]$. We then obtain that
\begin{eqnarray*}
&&\mbox{the first line of (\ref{eq:V2prebd})}
\\
&&\qquad \le\int_{[-\pi,\pi]^d}
\frac{\mathrm{d}^dk}{|k|^4}\int_{-\pi}^\pi\frac{\mathrm{d}\theta
}{2\pi}
\frac{O((1-{m/m_{\rm c}})^{-l+\epsilon})}{(|\theta
|+|k|^2)^{1
+\epsilon}}
\\
&&\qquad\le O\biggl(\biggl(1-\frac{m}{m_{\rm c}}\biggr)^{-l+\epsilon}\biggr)\int_{[-\pi
,\pi]^d}
\frac{\mathrm{d}^dk}{|k|^{4+2\epsilon}}=O\biggl(\biggl(1-\frac{m}{m_{\rm c}}\biggr)^{-l+
\epsilon}\biggr)
\end{eqnarray*}
if $\epsilon<\frac{d-4}2$. This completes the proof of (\ref
{eq:2<r<alpha-pibd})
for oriented percolation. This completes the proof of Lemma~\ref
{lmm:induction2}.
\end{pf*}

\subsection{\texorpdfstring{Proof of
Proposition~\protect\ref{proposition:r<alpha&2}}{Proof of Proposition~3.1}}
\label{ss:prop1}

First, we note that (\ref{eq:M-asy}) implies (\ref{eq:mu-asy}). To see
this, we
first substitute (\ref{eq:Jhat-diff}) and (\ref{eq:M-asy}) into (\ref
{eq:mu-def}) and
then use (\ref{eq:CICII}) [see (\ref{eq:M-source})] to obtain
\begin{eqnarray*}
\mu&=&\frac{m_{\rm c}\,\partial_m\hat J_{m_{\rm c}}(0) (1-
{m/m_{\rm c}})+O((1-{m/m_{\rm c}})^{1+\epsilon})}{M_{\mathrm c}v_\alpha+O
((1-
{m/m_{\rm c}})^\epsilon
)}
\\
&=&\frac{1-{m/m_{\rm c}}}{C_{\mathrm{II}}v_\alpha}+O
\biggl(\biggl(1-\frac{m}
{m_{\rm c}}\biggr)^{1+\epsilon}\biggr).
\end{eqnarray*}
Therefore, to complete the proof of Proposition~\ref{proposition:r<alpha&2},
it suffices to show (\ref{eq:M-asy}).

It is easier to prove (\ref{eq:M-asy}) for $\alpha\le2$. In this case,
$M$ in
(\ref{eq:M-def}) is reduced to
\begin{eqnarray*}
M=
\cases{
m,&\quad\mbox{RW/SAW},
\cr
\hat\pi^\mathrm{OP}(0)pm,&\quad\mbox{OP}.
}
\end{eqnarray*}
Therefore, (\ref{eq:M-asy}) is trivial for random walk and
self-avoiding walk.
For oriented percolation, we use (\ref{eq:critpt}) and (\ref{eq:property2})
to obtain
\begin{eqnarray*}
M_{\mathrm c}-M&=&\hat\pi_{\mathrm c}^\mathrm{OP}(0)p(m_{\rm
c}-m)+\bigl(\hat\pi_{\mathrm c}^\mathrm{OP}
(0)-\hat\pi^\mathrm{OP}(0)\bigr)
pm
\\
&=&1-\frac{m}{m_{\rm c}}+O(L^{-d}) \biggl(1-\frac{m}{m_{\rm c}}\biggr),
\end{eqnarray*}
where the $O(L^{-d})$ term is uniform in $m\le m_{\rm c}$.
This implies (\ref{eq:M-asy}).

It remains to prove (\ref{eq:M-asy}) for $\alpha>2$. In fact, we only need
investigate the crossover terms in (\ref{eq:M-def}) that are
proportional to
${\mathbh1}_{\{\alpha>2\}}$ and show that
%
\begin{eqnarray}\label{eq:desiredbd}
|\nabla_{1}^2\hat\pi_{\mathrm c}(0)-\nabla_{1}^2\hat\pi(0)|\le
O
\biggl(\biggl(1-\frac{m}
{m_{\rm c}}\biggr)^\epsilon\biggr)
\end{eqnarray}
since the above proof for $\alpha\le2$ directly applies to the noncrossover
terms. Notice that, for $\epsilon\in(0,1)$,
\begin{eqnarray*}
0&\le& m_{\rm c}^t-m^t\le m_{\rm c}^t\biggl(1-\biggl(\frac{m}{m_{\rm
c}}\biggr)^t\biggr)^{1-\epsilon}\biggl(
\frac{1-({m/m_{\rm c}})^t}{1-{m/m_{\rm c}}}
\biggr)^\epsilon\biggl(1-\frac{m}{m_{\rm c}}\biggr)^\epsilon
\\
&\le& m_{\rm c}^t t^\epsilon\biggl(1-\frac{m}{m_{\rm
c}}\biggr)^\epsilon
\end{eqnarray*}
so that
\begin{eqnarray*}
|\nabla_{1}^2\hat\pi_{\mathrm c}(0)-\nabla_{1}^2\hat\pi(0)|&\le&
\sum_{t\in
\mathbb{N}}(m_{\rm c}^t-m^t)
\sum_{x\in{\mathbb Z}^d}x_1^2 |\pi_t(x)|
\\
&\le&\biggl(1-\frac{m}{m_{\rm c}}\biggr)^\epsilon\sum_{t\in\mathbb
{N}}t^\epsilon
m_{\rm c}^t\sum_{x\in{\mathbb Z}^d}
x_1^2 |\pi_t(x)|.
\end{eqnarray*}
Moreover, since
\begin{eqnarray*}
t^\epsilon=\frac{t}{t^{1-\epsilon}}=\frac{t}{\Gamma(1-\epsilon
)}\int_0^\infty
\ell^{-\epsilon}e^{-\ell t}\, \mathrm{d}\ell,
\end{eqnarray*}
we have
\begin{eqnarray}\label{eq:desiredprebd}
&&|\nabla_{1}^2\hat\pi_{\mathrm c}(0)-\nabla_{1}^2\hat\pi(0)|\nonumber
\\[-8pt]\\[-8pt]
&&\qquad \le
\frac
{(1-{m/m_{\rm c}})
^\epsilon}{\Gamma(1-\epsilon)}\int_0^\infty\frac{\mathrm{d}\ell
}{\ell^\epsilon}
\sum_{t\in\mathbb{N}}t (m_{\rm c}e^{-\ell})^t\sum_{x\in{\mathbb
Z}^d}x_1^2|\pi_t(x)|.\nonumber
\end{eqnarray}
To show (\ref{eq:desiredbd}), it thus suffices to prove that the above
integral with respect to $\ell$ is $O(1)$ for sufficiently small
$\epsilon$.

First, we consider self-avoiding walk. By the diagrammatic bound on
$\pi_t^{(2)}(x)$ in (\ref{eq:saw-diagbds}) [see (\ref
{eq:pi2sawbd})], we
readily obtain
\begin{eqnarray*}
\sum_{t\in\mathbb{N}}t m^t\sum_{x\in{\mathbb Z}^d}x_1^2 \pi
_t^{(2)}(x)&\le&
\sum_{s,s',
s''\in\mathbb{N}}(s+s'+s'') m^{s+s'+s''}
\\
&&{}\qquad\hspace*{11pt}\times\sum_{x\in{\mathbb Z}^d}x_1^2 \varphi
_s^\mathrm{SAW}(x)
\varphi_{s'}^\mathrm{SAW}(x) \varphi_{s''}^\mathrm
{SAW}(x)
\\
&\le&3W^{(2)}\sum_{x\in{\mathbb Z}^d}\sum_{s,s'\in\mathbb{N}}s
m^{s+s'}\varphi_s^\mathrm{SAW}(x)
\varphi_{s'}^\mathrm{SAW}(x)
\\
&\le&3B'W^{(2)},
\end{eqnarray*}
where
%
\begin{eqnarray}\label{eq:B'-def}
B'\equiv B'(m)=\sup_{y\in{\mathbb Z}^d}\sum_{x\in{\mathbb Z}^d}\sum
_{t\in\mathbb{N}}t
m^t\varphi_t^\mathrm{SAW}(x)
\sum_{s=0}^\infty m^s\varphi_s^\mathrm{SAW}(y-x).
\end{eqnarray}
Similarly to the above and the derivation of (\ref{eq:piNsawbd}), we
can show
that, by (\ref{eq:saw-diagbds}) and~(\ref{eq:naive}),
\begin{eqnarray*}
\sum_{t\in\mathbb{N}}t m^t\sum_{x\in{\mathbb Z}^d}x_1^2 \pi
_t^{(N)}(x)\le N^4
\bigl(B^{(0)}\bigr)^{N-2}B'W^{(2)}\qquad[N\ge2].
\end{eqnarray*}
Since $B^{(0)}=O(L^{-d})$ and $W^{(2)}=O(1)$ uniformly in
$m\le m_{\rm c}$ if
$d>4$ [see formulas (\ref{eq:B0sawbd})--(\ref{eq:W2lsawbd})], we obtain that, for
$L\gg1$,
%
\begin{eqnarray}\label{eq:desiredprebd-saw1}
&&\int_0^\infty\frac{\mathrm{d}\ell}{\ell^\epsilon}\sum_{t\in
\mathbb{N}
}t (m_{\rm c}e^{-
\ell})^t\sum_{x\in{\mathbb Z}^d}x_1^2|\pi_t(x)|\nonumber
\\[-8pt]\\[-8pt]
&&\qquad \le\underbrace{\sum_{N=2}^\infty O(N^4)
O(L^{-d})^{N-2}}_{O(1)}\int_0^\infty
\frac{\mathrm{d}\ell}{\ell^\epsilon} B'(m_{\rm c}e^{-\ell
}).\nonumber
\end{eqnarray}

We now show that the integral of $B'(m_{\rm c}e^{-\ell})/\ell
^\epsilon$ is
uniformly bounded\vspace*{1pt} if $\epsilon<\frac{d-4}2$. First, we replace
$t\varphi_t^\mathrm{SAW}(x)$ in (\ref{eq:B'-def}) by the
following bound due to
subadditivity:
\begin{eqnarray*}
t\varphi_t^\mathrm{SAW}(x)=\sum_{s=1}^t\varphi
_t^\mathrm{SAW}(x)\le\sum
_{s=1}^t(\varphi_{s
-1}^\mathrm{SAW}*D*\varphi_{t-s}^\mathrm{SAW})(x).
\end{eqnarray*}
Then, by using $|\hat D(k)|\le1$ and (\ref{eq:2<r<alpha-varphibd})
for $n=0$,
we obtain
%
\begin{eqnarray}\label{eq:B'bd}
B'(m_{\rm c}e^{-\ell})&\le& m_{\rm c}e^{-\ell}\int_{[-\pi,\pi
]^d}\frac{\mathrm{d}^dk}{(2
\pi)^d} |\hat\varphi^\mathrm{SAW}(k,m_{\rm c}e^{-\ell})|^3\nonumber
\\[-8pt]\\[-8pt]
&\le& O(1)\int_{[-\pi,\pi]^d}\frac{\mathrm{d}^dk}{|k|^4} \frac
{e^{-\ell}}{1-
e^{-\ell}+|k|^2},\nonumber
\end{eqnarray}
where the $O(1)$ term is independent of $\ell$. However, for
$\epsilon\in(0,1)$,
\begin{eqnarray*}
\int_0^\infty\frac{\mathrm{d}\ell}{\ell^\epsilon} \frac{e^{-\ell
}}{1-e^{-\ell}
+|k|^2}&\le&\frac1{1-e^{-1}}\biggl(\int_0^1\frac{\mathrm{d}\ell
}{\ell^\epsilon}
\frac1{\ell+|k|^2}+\int_1^\infty\frac{\mathrm{d}\ell}{\ell
^\epsilon} e^{-\ell}\biggr)
\\
&\le&\frac1{1-e^{-1}}\biggl(\int_0^{|k|^2}\frac{\mathrm{d}\ell}{\ell
^\epsilon}
\frac1{|k|^2}+\int_{|k|^2}^1\frac{\mathrm{d}\ell}{\ell^{1+\epsilon
}}+1\biggr)
\\
&=&O(|k|^{-2\epsilon}).
\end{eqnarray*}
Therefore, if $\epsilon<\frac{d-4}2$, then we obtain
%
\begin{eqnarray}\label{eq:desiredprebd-saw2}
\int_0^\infty\frac{\mathrm{d}\ell}{\ell^\epsilon} B'(m_{\rm
c}e^{-\ell})\le O(1)
\int_{[-\pi,\pi]^d}\frac{\mathrm{d}^dk}{|k|^{4+2\epsilon}}=O(1).
\end{eqnarray}
Combining (\ref{eq:desiredprebd}), (\ref{eq:desiredprebd-saw1}) and
(\ref{eq:desiredprebd-saw2}), we complete the proof of (\ref
{eq:desiredbd}) for
self-avoiding walk.

For oriented percolation, similarly to the derivation of (\ref{eq:piNopbd}),
we can show that, for $N\ge0$,
\begin{eqnarray*}
&&\sum_{t\in\mathbb{N}}t m^t\sum_{x\in{\mathbb Z}^d}x_1^2 \pi
_t^{(N)}(x)
\\
&&\qquad\le
(N+1)^2\bigl(
T^{(0)}\bigr)^{N-2}\bigl(\bigl(N\bigl(1+T^{(0)}\bigr)+T^{(0)}
\bigr)
T^{(0)} V'
\\
&&{}\hspace*{88pt}\qquad\quad+N\bigl((N-1)\bigl(1+T^{(0)}\bigr)+3T^{(0)}\bigr)T' V^{(0)}\bigr),
\end{eqnarray*}
where
\begin{eqnarray*}
V'\equiv V'(m)
&=&\sup_{(x,t)\in{\mathbb Z}^{d+1}}\sum_{(y,s)\in{\mathbb
Z}^{d+1}}|y_1|^2(mD*
\varphi_s^\mathrm{OP})(y) m^s |s-t+1|\nonumber
\\
&&\hspace*{81pt}{}\times (D*\varphi
_{s-t}^\mathrm{OP})(y-x),\nonumber
\\
T'\equiv T'(m)
&=&\sup_{(x,t)\in{\mathbb Z}^{d+1}}\sum_{(y,s),(y',s')\in{\mathbb
Z}^{d+1}}(mD*
\varphi_s^\mathrm{OP})(y) m^s |s'-t+1|\nonumber
\\[-8pt]\\[-8pt]
&&\hspace*{105pt}{}\times  (D*\varphi
_{s'-t}^\mathrm{OP})(y'-x)\nonumber
\\
&&\hspace*{105pt}{}\times\bigl(\varphi_{s-s'}^\mathrm
{OP}(y-y')+\varphi_{s'-s}^\mathrm{OP}\nonumber
(y'-y)\bigr).
\end{eqnarray*}
Since $T^{(0)}=O(L^{-d})$ and $V^{(0)}|_{l=1}=O(1)$ uniformly in
$m\le m_{\rm c}$ if $d>4$ and $p\le p_{\mathrm c}$ [see formulas (\ref
{eq:T0bd}) and (\ref{eq:V0bd})], we obtain
that, for $L\gg1$,
\begin{eqnarray}\label{eq:desiredprebd-op1}
&&\int_0^\infty\frac{\mathrm{d}\ell}{\ell^\epsilon}\sum_{t\in
\mathbb{N}
}t (m_{\rm c}
e^{-\ell})^t\sum_{x\in{\mathbb Z}^d}x_1^2|\pi_t(x)|\nonumber
\\[-8pt]\\[-8pt]
&&\qquad \le O(1)\int_0^\infty\frac{\mathrm{d}\ell}{\ell^\epsilon}
\bigl(V'(m_{\rm c}
e^{-\ell})+T'(m_{\rm c}e^{-\ell})\bigr).\nonumber
\end{eqnarray}
However, by the Markov property,
\begin{eqnarray*}
(t+1)(D*\varphi_t^\mathrm{OP})(x)=\sum_{s=0}^t(D*\varphi
_t^\mathrm{OP})(x)\le\sum_{s=0}^t
(\varphi_s^\mathrm{OP}*D*\varphi_{t-s}^\mathrm{OP})(x).
\end{eqnarray*}
Applying this bound to the definitions of $V'$ and $T'$ and then using
$|\hat D(k)|\le1$ and (\ref{eq:2<r<alpha-varphibd}) for $n=0,2$, we obtain
\begin{eqnarray*}
\left.\begin{array}{c}
V'(m_{\rm c}e^{-\ell}) \\[2pt] T'(m_{\rm c}e^{-\ell})
\end{array}
\right\}\le O(1)\int_{[-\pi,\pi]^d}\frac{\mathrm{d}^dk}{|k|^4}
\frac{e^{-\ell}}
{1-e^{-\ell}+|k|^2}.
\end{eqnarray*}
Recalling (\ref{eq:B'bd}) and (\ref{eq:desiredprebd-saw2}), we conclude
that (\ref{eq:desiredprebd-op1})
is uniformly bounded. This completes the proof of (\ref{eq:desiredbd}) for
oriented percolation. We have thus completed the proof of
Proposition~\ref{proposition:r<alpha&2}.

\begin{appendix}
\section*{Appendix}

\subsection{Asymptotics of $1-\hat D(k)$ for small $k$}\label{appendix:D}

In this appendix, we will use the following notation for convenience:
\begin{eqnarray*}
|\hspace*{-1,5pt}\|x\|\hspace*{-1,5pt}|_\ell=|x|\vee\ell\qquad[\ell>0].
\end{eqnarray*}

\renewcommand{\thetheorem}{A.\arabic{theorem}}
\setcounter{theorem}{0}
\begin{lemma}\label{lemma:D}
Let $\alpha,\rho>0$ and
\begin{eqnarray*}
h(x)=\frac{1+O(|\hspace*{-1,5pt}\|x|\hspace*{-1,5pt}\|_1^{-\rho})}{|\hspace*{-1,5pt}\|x|\hspace*{-1,5pt}\|_1^{d+\alpha}}\qquad
[x\in{\mathbb R}^d].
\end{eqnarray*}
Suppose that $h$ is a rotation-invariant function. Then, there exist
$\epsilon>0$ and $v_\alpha=O(L^{\alpha\wedge2})$ such that, for $|k|<1/L$,
the 1-step distribution $D$ in (\ref{eq:kac}) satisfies
\renewcommand{\theequation}{A.\arabic{equation}}
\setcounter{equation}{0}
\begin{eqnarray}\label{eq:1-hatDasy}
1-\hat D(k)=v_\alpha|k|^{\alpha\wedge2}\times
\cases{
1+O((L|k|)^\epsilon),&\quad$\alpha\ne2$,
\cr
\log\dfrac1{L|k|}+O(1),&\quad$\alpha=2$.
}
\end{eqnarray}
\end{lemma}

\begin{pf}
The case for $\alpha>2$ is easy. By the Taylor expansion of
$1-\cos(k\cdot x)$ and using the ${\mathbb Z}^d$-symmetry of $D$,
\begin{eqnarray*}
1-\hat D(k)=\sum_{x\in{\mathbb Z}^d}\bigl(1-\cos(k\cdot x)\bigr)
D(x)=\frac
{|k|^2}{2d}
\sum_{x\in{\mathbb Z}^d}|x|^2D(x)+O((L|k|)^{2+\epsilon})
\end{eqnarray*}
holds provided that $0<\epsilon<2\wedge(\alpha-2)$. This proves
(\ref{eq:1-hatDasy})
with $v_\alpha\equiv\sigma^2/(2d)=O(L^2)$.

It remains to prove (\ref{eq:1-hatDasy}) for $\alpha\le2$. First, we
note that,
by definition,
\begin{eqnarray*}
D(x)=\frac{c_h}{L^d} h(x/L)\qquad[x\in{\mathbb Z}^d],
\end{eqnarray*}
where
\begin{eqnarray*}
c_h=\biggl(\frac1{L^d}\sum_{y\in{\mathbb Z}^d/L}h(y)
\biggr)^{-1}=\int_{{\mathbb R}^d
}h(y)\, \mathrm{d}^d
y+O(L^{-1}).
\end{eqnarray*}
Taking the Fourier transform yields
\begin{eqnarray*}
1-\hat D(k)&=&\frac{c_h}{L^d}\sum_{x\in{\mathbb Z}^d}\bigl(1-\cos
(k\cdot
x)\bigr) h\biggl(
\frac{x}L\biggr)
\\
&=&\frac{c_h}{(L|k|)^d}\biggl(|k|^d\sum_{y\in|k|{\mathbb Z}^d}
\bigl(1-\cos
(e_k\cdot y)\bigr)
 h\biggl(\frac{y}{L|k|}\biggr)\biggr),
\end{eqnarray*}
where $e_k=k/|k|$. By the Riemann sum approximation for small $k$ and the
rotational invariance of $h$, we obtain
\begin{eqnarray*}
1-\hat D(k)&=&\frac{c_h(1+O(|k|))}{(L|k|)^d}\int_{|y|\ge|k|}
\bigl(1-\cos(e_k
\cdot y)\bigr) h\biggl(\frac{y}{L|k|}\biggr)\, \mathrm{d}^dy
\\
&=&\frac{c_h(1+O(|k|))}{(L|k|)^d}\int_{|y|\ge|k|}(1-\cos y_1) h
\biggl(\frac{y}
{L|k|}\biggr)\, \mathrm{d}^dy
\\
&=&c_h(L|k|)^\alpha\bigl(1+O(|k|)\bigr)
\\
&&{}\times\int_{|y|\ge|k|}(1-\cos y_1)\biggl(\frac1{|\hspace*{-1,5pt}\|y|\hspace*{-1,5pt}\|_{L|k|}^{d+
\alpha}}+\frac{O((L|k|)^\rho)}{|\hspace*{-1,5pt}\|y|\hspace*{-1,5pt}\|_{L|k|}^{d+\alpha+\rho
}}\biggr)\,
\mathrm{d}^dy.
\end{eqnarray*}
This is the starting point of the analysis for $\alpha\le2$.

For $\alpha<2$, we note that
\begin{eqnarray*}
\int_{|y|\ge|k|}\frac{1-\cos y_1}{|\hspace*{-1,5pt}\|y|\hspace*{-1,5pt}\|_{L|k|}^{d+\alpha}}
\,\mathrm{d}^dy
&=&\int_{|y|\ge L|k|}\frac{1-\cos y_1}{|y|^{d+\alpha}}\, \mathrm{d}^dy
\\
&&{}+\underbrace{\int_{|k|\le|y|<L|k|}\frac{1-\cos
y_1}{(L|k|)^{d+\alpha}}
 \,\mathrm{d}^dy}_{O((L|k|)^{2-\alpha})}
 \\
&=&\int_{{\mathbb R}^d}\frac{1-\cos y_1}{|y|^{d+\alpha}} \,\mathrm{d}^dy-
\underbrace{\int_{|y|<L|k|}\frac{1-\cos y_1}{|y|^{d+\alpha}}\, \mathrm{d}^d
y}_{O((L|k|)^{2-\alpha})}
\\
&&{}+O((L|k|)^{2-\alpha}),
\end{eqnarray*}
where we have used $L|k|<1$ to estimate the error terms. Moreover,\vspace*{-2pt}
\begin{eqnarray*}
\int_{|y|\ge|k|}\frac{1-\cos y_1}{|\hspace*{-1,5pt}\|y|\hspace*{-1,5pt}\|_{L|k|}^{d+\alpha+\rho
}} \,\mathrm{d}^dy
&=&\underbrace{\int_{|y|\ge1}\frac{1-\cos y_1}{|y|^{d+\alpha+\rho
}}\, \mathrm{d}^d
y}_{O(1)}
+\int_{L|k|\le|y|<1}\frac{1-\cos y_1}{|y|^{d+\alpha+\rho}}
\,\mathrm{d}^dy
\\[-2pt]
&&{}+\underbrace{\int_{|k|\le|y|<L|k|}\frac{1-\cos
y_1}{(L|k|)^{d+\alpha+
\rho}} \,\mathrm{d}^dy}_{O((L|k|)^{2-\alpha-\rho})},\vspace*{-2pt}
\end{eqnarray*}
where\vspace*{-2pt}
\begin{eqnarray*}
\int_{L|k|\le|y|<1}\frac{1-\cos y_1}{|y|^{d+\alpha+\rho}}\, \mathrm{d}^dy=
\cases{
O(1),&\quad$\rho<2-\alpha$,\cr
O\biggl(\log\dfrac1{L|k|}\biggr),&\quad$\rho=2-\alpha$,\cr
O((L|k|)^{2-\alpha-\rho}),&\quad$\rho>2-\alpha$.
}\vspace*{-2pt}
\end{eqnarray*}
This proves (\ref{eq:1-hatDasy}) with $0<\epsilon<1\wedge(2-\alpha
)\wedge\rho$ and\vspace*{-2pt}
\begin{eqnarray*}
v_\alpha=c_hL^\alpha\int_{{\mathbb R}^d}\frac{1-\cos
y_1}{|y|^{d+\alpha}}
\,\mathrm{d}^dy.\vspace*{-2pt}
\end{eqnarray*}

For $\alpha=2$, we note that\vspace*{-2pt}
\begin{eqnarray*}
\int_{|y|\ge|k|}\frac{1-\cos y_1}{|\hspace*{-1,5pt}\|y|\hspace*{-1,5pt}\|_{L|k|}^{d+2}} \,\mathrm{d}^dy
&=&\underbrace{\int_{|y|\ge1}\frac{1-\cos y_1}{|y|^{d+2}}\, \mathrm
{d}^dy}_{O(1)}
+\int_{L|k|\le|y|<1}\frac{1-\cos y_1}{|y|^{d+2}} \,\mathrm{d}^dy
\\[-2pt]
&&{}+\underbrace{\int_{|k|\le|y|<L|k|}\frac{1-\cos
y_1}{(L|k|)^{d+2}}
\,\mathrm{d}^dy}_{O(1)}.\vspace*{-2pt}
\end{eqnarray*}
By the Taylor expansion of $1-\cos y_1$ and using $|y|^2=\sum_{j=1}^dy_j^2$,
we obtain\vspace*{-2pt}
\begin{eqnarray*}
\int_{L|k|\le|y|<1}\frac{1-\cos y_1}{|y|^{d+2}} \,\mathrm{d}^dy
&=&\frac12\int_{L|k|\le|y|<1}\frac{y_1^2}{|y|^{d+2}}\, \mathrm
{d}^dy+O(1
)\\[-2pt]
&=&\frac1{2d}\int_{L|k|\le|y|<1}\frac1{|y|^d}\, \mathrm{d}^dy+O(1)
\\[-2pt]
&=&\frac{\omega_d}
{2d}\log\frac1{L|k|}+O(1),\vspace*{-2pt}
\end{eqnarray*}
where $\omega_d\equiv2\pi^{d/2}/\Gamma(d/2)$ is the surface area of
the unit
$d$-sphere. Moreover,\vspace*{-2pt}
\begin{eqnarray*}
\int_{|y|\ge|k|}\frac{1-\cos y_1}{|\hspace*{-1,5pt}\|y|\hspace*{-1,5pt}\|_{L|k|}^{d+2+\rho}}
\,\mathrm{d}^dy
&=&\underbrace{\int_{|y|\ge L|k|}\frac{1-\cos y_1}{|y|^{d+2+\rho}}
\,\mathrm{d}^d
y}_{O((L|k|)^{-\rho})}
\\[-2pt]
&&{}+\underbrace{\int_{|k|\le|y|<L|k|}\frac{1-\cos
y_1}{(L|k|)^{d+2+\rho}}
\,\mathrm{d}^dy}_{O((L|k|)^{-\rho})}.\vspace*{-2pt}
\end{eqnarray*}
This proves (\ref{eq:1-hatDasy}) with $v_2=c_hL^2\omega_d/(2d)$.
\end{pf}

\subsection{Identity for the constant $K_r$}\label{appendix:K}
\begin{lemma}\label{lemma:K}
For $r\in(0,2)$,
%
\begin{eqnarray}\label{eq:lemma:K}
K_r\equiv\int_0^\infty\frac{1-\cos v}{v^{1+r}}\, \mathrm{d}v=\frac
\pi{2\Gamma(r+1)
\sin({r\pi}/2)}.
\end{eqnarray}
\end{lemma}

\begin{pf}
Below, we prove (\ref{eq:lemma:K}) only for $r\in(0,1]$. Since the
definition of
$K_r$ and the rightmost expression in (\ref{eq:lemma:K}) are both
analytic in
$r\in{\mathbb C}$ with \mbox{$0<\Re(r)<2$}, we can extend (\ref
{eq:lemma:K}) to $r\in(1,2)$
using analytic continuation.

First, we rewrite $K_r$ as
\begin{eqnarray}\label{eq:Kr-rewr}
K_r&=&\int_0^\infty\frac{\mathrm{d}u}{u^{1+r}}\int_0^u\sin v\, \mathrm{d}v
=\frac1r\int_0^\infty\frac{\sin v}{v^r}\,\mathrm{d}v\nonumber
\\[-8pt]\\[-8pt]
&=&\mathop{\lim_{R\to\infty}}_{\delta\to0}\frac1{2ir}\int
_\delta^R
\frac{e^{iv}-e^{-iv}}{v^r}\,\mathrm{d}v.\nonumber
\end{eqnarray}
For $a>0$, we let
\begin{eqnarray*}
\gamma_a^\pm&=&\biggl\{z=ae^{\pm i\theta}\dvtx \theta\mbox{ increases from 0
to }
\dfrac\pi2\biggr\},
\\
\eta^\pm&=&\{z=\pm iv\dvtx v\mbox{ increases from $\delta$ to }R\}.
\end{eqnarray*}
Then, by the Cauchy integral formula,
\begin{eqnarray*}
\int_\delta^R\frac{e^{iv}}{v^r}\,\mathrm{d}v&=&\int_{\gamma_\delta
^+}\frac{e^{iz}}
{z^r}\,\mathrm{d}z+\int_{\eta^+}\frac{e^{iz}}{z^r}\,\mathrm{d}z-\int
_{\gamma_R^+}
\frac{e^{iz}}{z^r}\,\mathrm{d}z
\\
&=&i\int_0^{\pi/2}\frac{e^{i\delta e^{i\theta}}}{(\delta
e^{i\theta})^{r-1}}
\,\mathrm{d}\theta+i^{1-r}\int_\delta^R\frac{e^{-v}}{v^r}\,\mathrm{d}v-
\underbrace{i\int_0^{\pi/2}\frac{e^{iRe^{i\theta
}}}{(Re^{i\theta})^{r-1}}\,
\mathrm{d}\theta}_{O(R^{-r})}.
\end{eqnarray*}
Similarly,
\begin{eqnarray*}
\int_\delta^R\frac{e^{-iv}}{v^r}\,\mathrm{d}v&=&\int_{\gamma_\delta^-}
\frac{e^{-iz}}{z^r}~\mathrm{d}z+\int_{\eta^-}\frac
{e^{-iz}}{z^r}\,\mathrm{d}z
-\int_{\gamma_R^-}\frac{e^{-iz}}{z^r}\,\mathrm{d}z
\\
&=&-i\int_0^{\pi/2}\frac{e^{-i\delta e^{-i\theta}}}{(\delta e^{-i
\theta})^{r-1}}\,\mathrm{d}\theta+(-i)^{1-r}\int_\delta^R\frac{e^{-v}}{v^r}\,
\mathrm{d}v+O(R^{-r}).
\end{eqnarray*}
Substituting these expressions back into (\ref{eq:Kr-rewr}) yields
%
\begin{eqnarray}\label{eq:Kr-rewr2}
K_r&=&\mathop{\lim_{R\to\infty}}_{ \delta\to0}\biggl(\frac{\delta
^{1-r}}{2r}
\int_0^{\pi/2}\biggl(\frac{e^{i\delta e^{i\theta}}}{e^{i\theta(r-1)}}+
\frac{e^{-i\delta e^{-i\theta}}}{e^{-i\theta(r-1)}}\biggr)\,\mathrm
{d}\theta\nonumber
\\[-9pt]\\[-9pt]
&&{}\hspace*{64pt}+i^{-r}\frac{1+(-1)^{-r}}{2r}\int_\delta^R\frac{e^{-v}}{v^r}\,\mathrm{d}v
\biggr).\nonumber
\end{eqnarray}
If $r=1$, then the second term is absent due to the cancelation $1+(-1)=0$.
By dominated convergence, we obtain
%
\begin{eqnarray}\label{eq:K1}
K_1=\lim_{\delta\to0}\frac12\int_0^{\pi/2}(e^{i\delta
e^{i\theta}}
+e^{-i\delta e^{-i\theta}})\,\mathrm{d}\theta=\int_0^{\pi/
2}\,\mathrm{d}
\theta=\frac\pi2.
\end{eqnarray}
If $r\in(0,1)$, on the other hand, the first term in (\ref
{eq:Kr-rewr2}) is
$O(\delta^{1-r})$ and therefore goes to zero as $\delta\to0$. Since
$(-1)^{-r}=(-1)^r=i^{2r}$ and $i^r+i^{-r}=2\cos\frac{r\pi}2$, we obtain
\begin{eqnarray*}
K_r=\frac{\cos({r\pi}/2)}r\int_0^\infty\frac
{e^{-v}}{v^r}\,\mathrm{d}v
=\frac{\cos({r\pi}/2)}r\Gamma(1-r).
\end{eqnarray*}
Using the well-known relations $\Gamma(1-r) \Gamma(r)=\pi/\sin
(r\pi)$ and
$r\Gamma(r)=\Gamma(r+1)$, we finally arrive at
\begin{eqnarray*}
K_r=\frac{\cos({r\pi}/2)}{r\Gamma(r)}~\frac\pi{\sin(r\pi)}
=\frac\pi{2\Gamma(r+1)\sin({r\pi}/2)}.
\end{eqnarray*}
This is also valid for $r=1$, due to (\ref{eq:K1}). This completes
the proof of Lem-\break ma~\ref{lemma:K}.
\end{pf}

\end{appendix}

\section*{Acknowledgments}
The first-named author is grateful to the staff of L-Station at
Hokkaido University for their
support and hospitality during a visit (July 23--August 5, 2008). The
second-named author is
grateful to Tai-Ping Liu and the Institute of Mathematics at Academia Sinica
in Taiwan for their support and hospitality during a visit (November
5--14, 2009). We would like to thank the referee for many valuable comments
regarding an earlier version of this paper.

%
%

\printaddresses

\end{document}